\documentclass{article}
\usepackage{amsfonts}
\usepackage{graphics}
\usepackage{amsmath}
\usepackage{amscd}
\usepackage{color}
\usepackage{euscript}
\usepackage{amssymb}

\usepackage{rlepsf}
\newtheorem{Theorem}{Theorem}

\newtheorem{Corollary}[Theorem]{Corollary}

\newtheorem{Example}[Theorem]{Example}

\newtheorem{Definition}[Theorem]{Definition}
\newtheorem{Remark}[Theorem]{Remark}

\newtheorem{Lemma}[Theorem]{Lemma}

\newtheorem{Proposition}[Theorem]{Proposition}

\newtheorem{claim}{Claim}

\newtheorem{Fundamental Theorem}{Fundamental Theorem}

\newtheorem{Claim}[claim]{Claim}

\newenvironment{Proof}[1][Proof]{\textbf{#1.} }{\ \rule{0.5em}{0.5em}}

\def \Wc {{\EuScript{W}}}
\def \U {\mathcal{U}}
\def \Der {\mathrm{Der}}
\def \op {\mathrm{op}}
\def \lH {{\mathfrak{H}}}
\def \fo {{\textrm{ for each }}}
\def \an  {{\textrm{ and }}}
\def \wh  {{\textrm{ where }}}
\def \de {\delta}
\def \ra {\xrightarrow}
\def \ad {\mathrm{ad}}
\def \ds {\frac{\d}{\d s}}
\def \dt {\frac{\d}{\d t}}

\def \dx {\frac{\d}{\d x}}
\def \du {\frac{\d}{\d u}}

\def \R {\mathbb{R}}
\def \ll {\mathfrak{l}}
\def \be {\natural}
\def \S {\mathcal{S}}
\def \Sc {\EuScript{S}}

\def \rank {\mathrm{Rank}}
\def \Dc {\mathcal{D}}
\def \lG {\mathfrak{G}}
\def \M {\EuScript{M}}
\def \A {\mathcal{A}}
\def \le {\mathfrak{e}}
\def \Gc {\mathcal{G}}
\def \Hc {\mathcal{H}}
\def \and {\textrm{ and }}
\def \g {\gamma}

\def \G {\Gamma}
\def \A {\mathcal{A}}

\def \C {\mathcal{C}}

\def \e {\epsilon}

\def \d {\partial}
\def \t {\triangleright}
\def \id {\mathrm{id}}
\def \a {\alpha}
\def \b {\beta}
\def \w {\omega}
\def \W {\Omega}

\def \f {\phi}
\def \i {\iota}

\def \lg {\mathfrak{g}}
\def \GL {\mathrm{GL}}
\def \gl {\mathfrak{gl}}
\def \H {{\bf\EuScript{H}}}
\begin{document}
\author{Jo\~{a}o  Faria Martins \\{ \small Departamento de Matem\'{a}tica }\\ {\small
Faculdade de Ci\^{e}ncias e Tecnologia,   Universidade Nova de Lisboa}\\
{\small Quinta da Torre,
2829-516 Caparica,
Portugal }\\\\ Roger Picken  \\ {\small Departamento de Matem\'{a}tica,} \\ {\small Instituto Superior T\'{e}cnico, TU Lisbon } \\ {\small Av. Rovisco Pais, 1049-001 Lisboa, Portugal} \\ {\it \small rpicken@math.ist.utl.pt}}

\title{The fundamental Gray 3-groupoid of a smooth manifold and local 3-dimensional holonomy based on a 2-crossed module}
\maketitle
\begin{abstract}
We define the thin fundamental Gray 3-groupoid $\EuScript{S}_3(M)$ of a smooth manifold $M$ and define (by using differential geometric data)  3-dimensional holonomies, to be smooth strict Gray 3-groupoid maps $\Sc_3(M) \to \C(\Hc)$, where $\Hc$ is a 2-crossed module of Lie groups and $\C(\Hc)$ is the Gray 3-groupoid naturally constructed from $\Hc$. As an application, we define Wilson 3-sphere observables.
\end{abstract}

\noindent{\bf 2000 Mathematics Subject Classification:} {\it 53C29 
(primary); %Issues of holonomy,
18D05, % Double categories, $2$-categories, bicategories and generalisations}  
70S15 %Yang-Mills and other gauge theories 
(secondary) 
}
\section*{Introduction}
This is the first of a series of papers on (Gray) 3-bundles and 3-dimensional holonomy (also called parallel transport), aimed at categorifying the notion of 2-bundles and non-abelian gerbes with connection, and their 2-dimensional parallel transport; see \cite{BS,SW3,BrMe,FMP2,H,MP,ACG,Wo,Bar}. 
The main purpose {here} is to clarify the notion of 3-dimensional holonomy based on a Lie 2-crossed module, extending some of the constructions in  \cite{SW1,SW2,SW3,FMP1} to the fundamental Gray 3-groupoid of a smooth manifold (which we will construct in this article). For some related work see \cite{SSS,J1,J2,Br}.

The definition of a Gray 3-groupoid appears, for example, in \cite{Cr,KP,GPS}. The first  main result of this article concerns the construction of the fundamental (semi-strict) Gray 3-groupoid $\Sc_3(M)$ of a smooth manifold $M$, which is not obvious. 
The main innovation for defining  $\Sc_3(M)$ rests on the notion of laminated rank-2 homotopy, a weakening of the notion of rank-2 homotopy (see \cite{MP,FMP1,SW3,BS}), which makes 2-dimensional holonomy based on a pre-crossed module, as opposed to a crossed module, invariant. This also permits us to define the fundamental pre-crossed module of a smooth manifold $M$. 

The definition of a 2-crossed module is due to Conduch\'{e}; see \cite{Co} and \cite{KP,MuPo,Po}.  It is a complex of (not necessarily abelian) groups
$$L \ra{\de} E\ra{\d} G$$ 
together with  left actions  $\t$ by automorphisms of $G$ on $L$ and  $E$,  and a $G$ equivariant function {$\left\{,\right \} \colon E \times E \to L$} (called the Peiffer lifting), satisfying certain properties (Definition \ref{2cmlg}). The extension to 2-crossed modules of Lie groups is the obvious one. {{It is well known that  Gray 3-groupoids are modelled by 2-crossed modules (of groupoids); see \cite{KP}.}  We provide a detailed description of this connection in subsection \ref{nnn}.

In the light of this connection, given a Lie 2-crossed module $\H$, it is natural to define a (local) 3-dimensional holonomy as being a smooth (strict) Gray 3-functor $\Sc_3(M) \to \C(\H)$, where $\C(\H)$ is the  Gray 3-groupoid, with a single object, constructed out of $\H$.

{A differential 2-crossed module, also called a 2-crossed module of Lie algebras, see \cite{AA,E},} is given by a complex of Lie algebras: $$\ll \ra{\de} \le\ra{\d} \lg$$
together with  left actions $\t$ by derivations of $\lg$ on $\ll$ and $\le$,  and a $\lg$-equivariant {bilinear map} $\left\{,\right \} \colon \le \times \le \to \ll$
(called the Peiffer lifting), satisfying appropriate conditions (Definition \ref{d2cm}). Any 2-crossed module of Lie groups $\H$ defines, in the natural way, a differential 2-crossed module $\mathfrak{H}$, and this assignment is functorial.

The second main result of this paper concerns how to define 3-dimensional holonomies  $\Sc_3(M) \to \C(\H)$ by picking Lie-algebra valued differential forms $\w \in \A^1(M,\lg)$, $m \in \A^2(M,\le)$ and $\theta \in \A^3(M,\ll)$, satisfying  $\d(m)=d \w+[\w,\w]\doteq\W$, the curvature of $\w$, and $\de(\theta)=d m+\w \wedge^\t m\doteq \M$ (see the appendix for the  $\wedge^\t$ notation), the exterior covariant derivative of $m$, called the 2-curvature $3$-form of $(\w,m)$. This is analogous to the construction in  \cite{BS,FMP1,SW2}. A lot of the proofs will make use of Chen integrals in the loop space \cite{Ch}, {which is the approach taken in \cite{BS}.}

The 3-curvature 4-form $\Theta$ of the triple $(\w,m,\theta)$, satisfying the above equations, is defined as $\Theta=d \theta+\w\wedge^\t \theta-m \wedge^{\{,\}} m$ (see the appendix for the  $ \wedge^{\{,\}} $ notation). We will prove a relation between 3-curvature and 3-dimensional holonomy, completely analogous to that for principal $G$-bundles and 2-bundles {with a structural} crossed module of groups, including an Ambrose-Singer type theorem for {the triple $(\w,m,\theta)$}. This will also prove that the 3-dimensional holonomy is invariant under rank-3 homotopy, as long as it restricts to a  laminated rank-2 homotopy on the boundary (see Definition \ref{3strong}).

As in the case of 2-bundles, the 1-, 2- and 3-gauge transformations would be better understood by passing to the notion of a Gray triple groupoid, which (to keep the size of this paper within limits) we will analyse in a future article, where we will also  address the general definition of a Gray 3-bundle, and describe the corresponding $3$-dimensional parallel transport, categorifying the results of \cite{SW3,FMP2,MP,P}.

The 3-dimensional holonomy which we define in this article can be associated to embedded oriented 3-spheres $S$ in a manifold $M$, yielding a Wilson 3-sphere observable $\Wc(S,\w,m,\theta)\in \ker \de \subset L$ independent of the parametrisation of $S$ chosen, up to acting by elements of $G$. This will be a {corollary} of the invariance of the 3-dimensional holonomy under rank-3 homotopy, with laminated boundary.

 {Considering the Lie 2-crossed module given by a finite type chain complex of  vector spaces (see \cite{KP} and \ref{l3cc}), the construction in this article will describe (locally and as a Gray 3-functor) the first three instances of the holonomy of a representation up to homotopy; see \cite{AC}. Notice however that the  construction  in this article is valid for any Lie group 2-crossed module. We expect that to describe all instances of the $\w$-parallel transport of a representation up to homotopy one will need  Gray $\w$-groupoids and  2-crossed complexes. We hope to address this in a future paper.}

\tableofcontents
\section{2-crossed modules (of Lie groups and Lie algebras)}
\subsection{Pre-crossed modules (of Lie groups and Lie algebras)}
\subsubsection{Pre-crossed modules of Lie groups}
\begin{Definition}[Lie pre-crossed module]
A (Lie group) pre-crossed module $\Gc=(\d\colon E \to G, \t) $ is given by a Lie group map $\d\colon E \to G$ together with a smooth left action $\t$ of $G$ on  $E$ by automorphisms such that:
$${\d(g \t e)=g \d(e)g^{-1} , \textrm{ for each } g \in G \textrm{ and }e \in E.}$$
The Peiffer commutators in a pre-crossed module are defined as $$\left <e,f\right >=efe^{-1}\left(\d(e) \t f^{-1} \right) , \wh e,f \in E.$$
A pre-crossed module is said to be a crossed module if all of its Peiffer commutators {are trivial}, which is to say that:
$$ \d(e) \t f=efe^{-1} \textrm{ for each } e,f \in E.$$
\end{Definition}
Note that the map $(e,f) \in E \times E \mapsto \left <e,f \right > \in E$, {called the Peiffer pairing,} is $G$-equivariant:
$${g \t \left <e,f \right >= \left <g \t e, g \t f\right > , \fo e,f \in E \an g \in G.}$$
Moreover {$\left <e,f\right >=1$}  if either $e$ or $f$ is $1$. Therefore, the second differential of the Peiffer pairing defines a bilinear map $\left <,\right >\colon\le \times \le \to \le$, {where $\le=T_1 E$ is the Lie algebra of $E$.}

\subsubsection{Differential pre-crossed modules}
The infinitesimal counterpart of a Lie pre-crossed module is a differential pre-crossed module, also called a \emph{pre-crossed module of Lie algebras.}
\begin{Definition}[Differential pre-crossed module]
A differential pre-crossed module $\Gc=(\d\colon \le\to \lg,\t)$ is given by a Lie algebra map $\d\colon \le \to \lg$ together with a left action $\t$ of $\lg$ on $\le$ by derivations such that:
\begin{enumerate}
 \item $\d(X \t v)=[X,\d(v)]$, for each $v \in \le$ and $X \in \lg$.
\end{enumerate}
Note that the map  {$(X,v) \in \lg \times \le \mapsto X \t v \in \le$ is necessarily bilinear.} In addition we have $$X\t[v,w]=[X \t v,w]+[v,X \t w], \fo  v,w \in \le \an  X \in \lg$$ (this expresses the condition that $\lg$ acts on $\le$ on the left by derivations.)
\end{Definition}
Therefore if $\Gc$ is a Lie pre-crossed module {then} the induced structure on Lie algebras is a differential pre-crossed module $\lG$. Reciprocally, if $\lG$ is a differential pre-crossed module there exists a unique Lie pre-crossed module of simply connected Lie groups $\Gc$ whose differential form is $\lG$ (any action by derivations of a Lie algebra on another Lie algebra can always be lifted to an action by automorphisms, considering the corresponding simply connected Lie groups).

In a differential pre-crossed module, the Peiffer commutators are defined as: 
$$\left <u,v \right >=[u,v]-\d(u) \t v \in \le; \quad u, v \in \le.$$
The map $ (u,v) \in \le \times \le \mapsto \left <u,v\right > \in \le$ (the Peiffer pairing) is bilinear  (though not necessarily symmetric),  and  coincides with the second differential of the Peiffer pairing in $E$.
In addition the Peiffer pairing is $\lg$-equivariant:
$${X \t \left <u,v\right >=\left <X \t u, v\right >+\left <u,X\t v\right >, \fo X \in \lg, \an u,v \in \le .}$$
A differential pre-crossed module is said to be a differential crossed module if all of its Peiffer commutators vanish, which is to say that:
$$\d(u) \t v= [u,v]; \quad \fo u, v \in \le. $$
For more on crossed modules of groups,  Lie groups and Lie algebras, see \cite{B,BL,BC,BH1,BHS,FM}.

\subsection{Lie 2-crossed modules and Gray 3-groupoids}\label{2X}
\subsubsection{2-Crossed modules of Lie groups}
We follow the conventions of \cite{Co} for the definition of a 2-crossed module. See also \cite{MuPo,KP,BG,Po,RS}.
\begin{Definition}[2-crossed module of Lie-groups]\label{2cmlg}
A 2-crossed module (of Lie groups) is given by a complex of Lie groups:
$$L \ra{\de} E\ra{\d} G$$
together with  smooth left actions $\t$ by automorphisms of $G$ on $L$ and  $E$ (and on $G$ by conjugation), and a $G$-equivariant smooth function {$\left\{,\right \} \colon E \times E \to L$} (called the Peiffer lifting). {As before $G$-equivariance means}
\begin{equation*}
a \t \{e,f\}= \{a \t e, a \t f\}, \fo a \in G \an e,f \in E.
\end{equation*}
{These are to satisfy:}
\begin{enumerate}
 \item $L \ra{\de} E\ra{\d} G$ is a complex of $G$-modules (in other words $\d$ and $\de$ are $G$-equivariant {and $\d \circ \de =1$}.)
\item $\de(\left\{e,f\right \})=\left <e,f\right >,$ for each $e,f \in E$. {Recall $\left <e,f\right >=efe^{-1} \d(e) \t f^{-1}. $}
\item $[l,k]=\left\{\de(l),\de(k)\right \}, $ for each $l,k \in L$. Here $[l,k]=lkl^{-1}k^{-1}$.
\item\label{B} $\left\{ef,g\right \}=\left\{e,fgf^{-1}\right \}\d(e)\t \left\{f,g\right \}$, for each $e,f,g \in E$.

\item \label{C}$\left\{e,fg\right \}=\left \{e,f\right\}\left \{e,g\right\}\left \{{\left <e,g\right >}^{-1}, \d(e) \t f\right\}$, where {$e,f,g \in E$}.
\item $\left\{\de(l),e\right \}\left\{e, \de(l)\right \}=l(\d(e) \t l^{-1})$, for each $e\in E$ and $l \in L$.
\end{enumerate}

\end{Definition}
Define:
\begin{equation}\label{tprime}
e \t'l=l\left\{\de(l)^{-1},e\right\}, \wh l\in L \an e \in E.
\end{equation}
 It follows from the previous axioms that $\t'$ is a left action of $E$ on $L$ by automorphisms (this is not entirely immediate; a proof is in \cite{Co,BG}). {There it is also shown that,} together with the map $\de\colon L \to M$, this defines a crossed module. This will be of prime importance later.  For the time being, note that condition 5. yields:
\begin{align}
 \{e,fg\}&=\{e,f\} (\d(e) \t f) \t'\{e,g\}=\left(\de(\{e,f\}) \t' (\d(e) \t f) \t'\{e,g\}\right) \{e,f\}\nonumber \\&=\left((efe^{-1}) \t' \{e,g\} \right)\{e,f\}.\label{ky}
\end{align}
{This proves that the definition of 2-crossed modules appearing here is equivalent to the one of \cite{KP,Po}.}
{The following simple lemma is useful here and later on:}
\begin{Lemma}\label{trivial}
In a 2-crossed module  {$\Hc=(L\ra{\de} E \ra{\d} G, \t, \{,\})$ we have $\left\{1_E,e\right \}=\left\{e,1_E\right \}=1_L$, for each $e \in E$.}
\end{Lemma}
\begin{Proof}
{Apply axioms \ref{B} and \ref{C} of the definition of a 2-crossed module  to {$\left\{1_E1_E,e\right \}$ and $\left\{e,1_E1_E\right \}$.}}
\end{Proof}

{Using this together with equations 4. and 5. of the definition of a 2-crossed module as well as equation (\ref{ky})  it follows (compare also with equation (\ref{invv})):}
\begin{Lemma}\label{Inverses}
For each $e,f \in E$ we have:
$$\{e,f\}^{-1}=\d(e) \t \{e^{-1},efe^{-1}\}, $$
$$\{e,f\}^{-1}=(efe^{-1}) \t'\{e,f^{-1}\}, $$
$$\{e,f\}^{-1}=(\d(e) \t f) \t' \{e,f^{-1}\}. $$
\end{Lemma}

By using condition 6. of the definition of a 2-crossed module it follows that:
\begin{equation}\label{jkl}
\d(e) \t l=(e \t' l) \{e,\de(l)^{-1}\}.
\end{equation}
Note:
\begin{align*} \big(\d(e) \t f\big) \t' \big(\d(e) \t l\big)&\doteq (\d(e) \t l) \{\d(e) \t \de(l^{-1}), \d(e) \t f\}\\&=( \d(e) \t l) \d(e) \t \{  \de(l^{-1}), f\}\\
&=\d(e) \t (f \t' l),
\end{align*}
{where we have used the fact that the Peiffer lifting  $\{,\}$ is $G$-equivariant and that $G$ acts on $L$ by automorphisms.}
We thus have the following identity for each $e,f \in E$ and $l \in L$:
\begin{equation}\label{wso}
  \big(\d(e) \t f\big) \t' \big(\d(e) \t l\big)=\d(e) \t (f \t' l).
\end{equation}

We also have:
\begin{align*}
(e \t' \{f,g\} ) \{e,\d(f) \t g\}&=\d(e) \t \{f,g\}\{e,(\d(f) \t g)fg^{-1}f^{-1}\}^{-1}\{e,\d(f) \t g\}\\
&=\d(e) \t \{f,g\} (\d(e) \t \d(f) \t g)\t' \{e,fg^{-1}f^{-1}\}^{-1}\\
&=\d(e) \t \{f,g\} \big(\d(e) \t \d(f) \t g) \t' (\d(e) \t( fg^{-1}f^{-1})\big)\t' \{e,fgf^{-1}\}\\
&=\d(e) \t \{f,g\} \big(\d(e) \t \de\{f,g\}^{-1}\big)\t' \{e,fgf^{-1}\}\\
&=\{e,fgf^{-1}\}\d(e) \t \{f,g\}\\
&=\{ef,g\},
\end{align*}
using equation (\ref{jkl})  in the first step, condition 5. of the definition of a 2-crossed module in the second, the third equation of Lemma \ref{Inverses} in the third. The penultimate step follows from the fact that $(\de \colon L \to E ,\t')$ is a crossed module.

{Let us display all equations we have proved (note that equation (\ref{imppp}) appears in \cite[page 162]{Co}):}
\begin{Lemma}
 In a 2-crossed module we have, for each $e,f,g \in E$:
\begin{equation}\label{imppp}
\{ef,g\}= (e \t' \{f,g\} ) \{e,\d(f) \t g\},  \end{equation}
\begin{equation}\label{condunche}
\{ef,g\}= \left\{e,fgf^{-1}\right \}\d(e)\t \left\{f,g\right \}, \end{equation}
and
\begin{equation}\{e,fg\}=\{e,f\} (\d(e) \t f) \t'\{e,g\},\end{equation}
 \begin{equation}\label{Porter}
\{e,fg\}=\left((efe^{-1}) \t' \{e,g\} \right)\{e,f\}.
\end{equation}
\end{Lemma}
In particular, by using the first equation, we have (compare with {Lemma} \ref{Inverses}):
\begin{equation}\label{invv}
\{e,f\}^{-1}=e \t' \{e^{-1},\d(e) \t f\}.
\end{equation}

Consider the totally intransitive groupoid with morphisms  $G \times L$ and objects $G$, the source and target maps being given by $(g,l) \mapsto g$, and identity as $g \mapsto (g,1_L)$. As composition,  we take the group multiplication in $L$. Consider also the groupoid with objects $G$ and morphisms $G \times E$, and source and target given by $(g,e) \mapsto g$ and $(g,e) \mapsto \d(e)^{-1}g$, respectively. The composition is 
$$ (g \ra{(g,e)} \d(e)^{-1} g \ra{(\d(e)^{-1}g,f)} \d(f)^{-1}\d(e)^{-1}g) =(g \ra{(g,ef)} \d(ef)^{-1}g).
$$ The map $\de\colon G \times L \to G \times E$ defined as $\de(g,k)=(g,\de(k))$ is a groupoid map and together with the left action by automorphisms of $G\times E$ on $G \times L$: $$ (g,e) \t' (\d(e)^{-1} g,l)=(g,e\t' l), \wh g \in G,  e \in E \an l \in L$$
defines a crossed module of groupoids; see \cite{B1,BH1,BHS,FMPo,No}{, a part of what is called a {\it braided regular crossed module} in \cite{BG}.}

\begin{Example}\label{trivialexample}
The simplest non-trivial example of a Lie 2-crossed module is probably the following one. Let $G$ be a Lie group. Consider
$$G \ra{\de} G \rtimes^{\ad} G \ra{\d} G, $$
where the multiplication in the semidirect product is: $ (e,f)(g,h)=(efgf^{-1},fh)$ and also $\d(g,h)=gh$. In addition put $\de(g)=(g^{-1},g)$. Here $e,f,g,h \in G$. The action of $G$ on $G$ is the adjoint action and on $G \rtimes^{\ad} G $ is $g \t (a,b)=(gag^{-1},gbg^{-1})$. The Peiffer lifting is: 
$$\left \{ (a,b),(c,d)\right\} =  [bdb^{-1},a], $$
where  $[x,y]=xyx^{-1}y^{-1}$.
Therefore $$\d(a,b) \t g=(ab) g (ab)^{-1} \an (a,b) \t' g=bgb^{-1}.$$ 
\end{Example}
\subsubsection{Differential 2-crossed modules}
{To transport 2-crossed modules of Lie groups to the Lie algebras world we use Lemma \ref{trivial}. This tells us that the second differential of the Peiffer lifting {of $\Hc$}  defines a bilinear map} $$\left\{,\right \}\colon\le \times \le \to \ll,$$
which is $\lg$-equivariant:
$$X \t \left\{u,v\right \}=\left\{X \t u, v\right \}+\left\{u,X\t v\right \}, \fo X \in \lg \an \textrm{each } u,v \in \le .$$
{Further relations motivated by the remaining properties of the Peiffer lifting hold. These can be gathered inside the definition of a differential  2-crossed module (also called a \emph{2-crossed module of Lie algebras}).  This definition appeared in \cite{E}. Note that our conventions are different.}
\begin{Definition}[Differential 2-crossed module]\label{d2cm}
A differential 2-crossed module is given by a  complex of Lie algebras:
$$\ll \ra{\de} \le\ra{\d} \lg$$
together with  left actions $\t$ by derivations of $\lg$ on $\ll$, $\le$ and $\lg$ ({on the latter} by the adjoint representation), and a {$\lg$-equivariant bilinear map} $\left\{,\right \} \colon \le \times \le \to \ll$:
$$X \t \left\{u,v\right \}=\left\{X \t u, v\right \}+\left\{u,X\t v\right \}, \fo X \in \lg \an u,v \in \le, $$
(called the Peiffer lifting) such that:
\begin{enumerate}
 \item $L \ra{\de} E\ra{\d} G$ is a complex of $\lg$-modules.
\item $\de(\left\{u,v\right \})=\left <u,v\right >,$ for each $u,v \in \le$. {(Recall $\left <u,v\right >=[u,v]-\d(u) \t v$.)}
\item $[x,y]=\left\{\de(x),\de(y)\right \}, $ for each $x,y \in \ll$. 
\item $\left\{[u,v],w\right \}=\d(u)\t\left\{v,w\right \}+\left\{u,[v,w]\right \}-\d(v)\t\left\{u,w\right \}-\left\{v,[u,w]\right \}$, for each $u,v,w \in \le$.
This is the same as: 
\begin{equation}\label{sssubs} 
\left\{[u,v],w\right \}=\{ \d(u) \t v,\w\}-\{ \d(v) \t u,\w\}-\{u,\de\{v,w\}\}+\{v,\de\{u,w\}\}
\end{equation}
\item $\left\{u,[v,w]\right \}=\left\{\de\left\{u,v\right \},w\right \}-\left\{\de\left\{u,w\right \},v\right \}$, for each $u,v,w \in \le$. This implies that {$\t'$ defined by }$v\t'x=-\left\{\de(x),v\right \}$ is a left action of $\le$ on $\ll$; see below.
\item $\left\{\de(x),v\right \}+\left\{v,\de(x)\right \}=-\d(v) \t x$, for each $x \in \ll$ and $v \in \le$.
\end{enumerate}
 \end{Definition}
{Analogously to the 2-crossed module of Lie groups  case we have:}
\begin{Lemma}
The action $\t'$ of $\le$ on $\ll$ is by derivations, and together with the map $\de\colon \ll \to \le$ defines a differential crossed module.
\end{Lemma}
We divide the proof into four claims.
\begin{Claim}
$(u,x) \in \le \times \ll \mapsto u \t' x =-\{\de(x),u\}$ is a left  action of $\le$ on $\ll$.
\end{Claim}
\begin{Proof}
For each $u,v \in \le$ and $x\in \ll$ we have:
$$[u,v] \t' x=-\{\de(x),[u,v]\}= \{\de\{\de(x),v\},u\}-\{\de\{\de(x),u\},v\}=u\t'(v\t'x)-v\t'(u\t'x)$$
\end{Proof}
\begin{Claim}
The left action $\t'$ of $\le$ on $\ll$ is by derivations.
\end{Claim}
\begin{Proof}
We need to prove that $u\t'[x,y]=[u\t'x,y]+[x,u \t' y]$, for each $u\in \le$ and $x,y \in \lg$. We have:
$$u\t'[x,y]=-\{\de([x,y]),u\}=-\{[\de(x),\de(y)],u\}=\{\de(y),[\de(x),u]\}-\{\de(x),[\de(y),u]\}.$$
We have used $4.$ together with the fact $\d\de=0$. The last term can be simplified as:
\begin{align*}
\{\de(y),[\de(x),u]\}-\{\de(x),[\de(y),u]\}&=\{\de(\{\de(y),\de(x)\}), u\})-\{\de(\{\de(y),u\}), \de(x)\})\\&-\{\de(\{\de(x),\de(y)\}), u\})+\{\de(\{\de(x),u\}), \de(y)\})\\&=-[x,u\t' y]-[u\t' x,y]+2u \t'\{\de(x),\de(y)\}\\&=-[x,u\t' y]-[u\t' x,y]+2u \t'[x,y].
\end{align*}
Therefore
$$ u\t'[x,y]=-[x,u\t' y]-[u\t' x,y]+2u \t'[x,y],$$
from which the result follows.
\end{Proof}
\begin{Claim}
$\de(x) \t' y=[x,y]$, for each $x,y \in \ll$.
\end{Claim}
\begin{Proof}
$$\de(x) \t' y=-\{\de(y),\de(x)\}=[x,y], \fo x,y \in \ll$$
\end{Proof}
\begin{Claim}
$\de( u \t' x)=[u,\de(x)]$, for each $x \in \ll$ and $u \in \le$.
\end{Claim}
\begin{Proof}
$$\de( u \t' x)=-\de(\{\de(x),u\})=-<\de(x),u>=[u,\de(x)],$$
since $\d \de=0$.
\end{Proof}
\subsubsection{{Lie 2-crossed modules and differential 2-crossed modules}}

{The definition of a differential 2-crossed module is an exact differential replica of the definition of a 2-crossed module of Lie groups. Straightforward calculations prove that:}
\begin{Theorem}\label{Lie}
Let $\Hc=(L \ra{\de} E\ra{\d} G,\t,\{,\})$ be a 2-crossed module of Lie groups. The induced chain complex of Lie algebras $\ll \ra{\de} \le\ra{\d} \lg$, together with the induced actions  $\t$ of $\lg$ on $\le$ and $\ll$  and the second differential of the Peiffer {lifting} $\{,\}$ defines a differential 2-crossed module $\lH$. Moreover the assignment $\Hc\mapsto \lH$ is functorial.
\end{Theorem}
{Theorem \ref{Lie} can be proved by taking differentials, in the obvious way, of the equations appearing in the  definition of a  2-crossed module of Lie groups. There is, however, a more conceptual way to prove this theorem, which {guarantees} that we can go in the opposite direction. Namely,  it is well known that the categories of simplicial groups with Moore complex of length {two} and of 2-crossed modules are equivalent, see for example \cite{Co,P,BG}. This equivalence of categories also holds in the Lie algebra case, {as is proved} in \cite{E}. From standard Lie theory it follows that the categories of simplicial Lie algebras and simplicial (simply connected) Lie groups are equivalent. Taking Moore complexes, it therefore follows that:}

\begin{Theorem}
{Given a differential 2-crossed module $\lH=(\ll \ra{\de} \le\ra{\d} \lg,\t,\{,\})$  there exists a 2-crossed module of  Lie groups $\Hc$ whose differential form is $\lH$. }
\end{Theorem}
\subsubsection{The definition of Gray 3-groupoids}\label{dg3g}
We now define Gray 3-groupoids. Our conventions are slightly different from the ones of \cite{KP,Cr}.

A (small) Gray 3-groupoid $\C$ is given by a set $C_0$ of objects, a set $C_1$ of morphisms, a set $C_2$ of 2-morphisms and a set $C_3$ of $3$-morphisms, and maps $\d^\pm_i \colon C_k \to C_{i-1}$, where $i=1,\ldots, k$ {(and $k=1,2,3$)} such that:
\begin{enumerate}

\item $\d_2^\pm\circ \d_3^\pm=\d_2^\pm$, as maps $C_3 \to C_1$.
\item $\d_1^\pm=\d_1^\pm\circ \d_2^\pm=\d_1^\pm \circ \d_3^\pm,$  as maps $C_3 \to C_0$.
\item $\d_1^\pm=\d_1^\pm\circ \d_2^\pm$, as maps $C_2 \to C_0$. 
 \item There exists an upwards multiplication $J \be_3 J'$ of  3-morphisms if $\d_3^+(J)=\d_3^-(J')$, making $C_3$ into a groupoid whose set of objects is $C_2$ (identities are implicit).
\item There exists a vertical composition  $$\G \be_2 \G'=\begin{pmatrix} \G' \\ \G\end{pmatrix} $$ of  2-morphisms if $\d_2^+(\G)=\d_2^-(\G')$, making $C_2$ into a groupoid whose set of objects is $C_1$ (identities are implicit).

\item There exists a vertical composition $$J \be_2 J'=\begin{pmatrix} J'\\ J\end{pmatrix} $$ of 3-morphisms whenever $\d_2^+ (J)=\d_2^- (J')$ making the set of 3-morphisms into a groupoid with set of objects $C_1$ and such that the boundaries $\d_3^\pm \colon C_3 \to C_2 $ are functors.
\item The vertical and upwards compositions of 3-morphisms satisfy the interchange law $(J\be_3 J')\be_2 (J_1\be_3 J'_1)=(J\be_2 J_1)\be_3 (J'\be_2 J'_1)$, {whenever the compositions are well defined}. Combining with the previous axioms, this means that the vertical and upwards compositions of 3-morphisms and the vertical composition of $2$-morphisms  give $C_3$ the structure of a 2-groupoid, with set of objects being $C_1$, set of morphisms $C_2$ and set of 2-morphisms $C_3$. (The definition of a 2-groupoid appears for example in \cite{HKK}. It is well known that the categories of (small) 2-groupoids and of crossed modules of groupoids are equivalent; see for instance \cite{BHS,BS}.)

\item {\bf (Existence of whiskering by 1-morphisms)} For each $x,y$ in $C_0$ we can therefore define a 2-groupoid $\C(x,y)$ of all $1$-, 2- and 3-morphisms $b$ such that {$\d_1^-(b)=x$} and $\d_1^+(b)=y$. Given a 1-morphism $\g$ with $\d^-_1(\g)=y$ and $\d^+_1(\g)=z$ there exists a 2-groupoid map {$\be_1 \g\colon \C(x,y) \to \C(y,z)$}, called right whiskering. Similarly if $\d^+_1(\g')=x$ and $\d^-_1(\g')=w$ there exists a 2-groupoid map $\g' \be_1\colon \C(x,y) \to \C(w,y),$ called left whiskering.

\item There exists therefore a horizontal composition of $\g\be_1 \g'$ of 1-morphisms if $\d_1^+(\g)=\d_1^-(\g')$, which is to be associative and to define a groupoid with set of objects $C_0$ and set of morphisms $C_1$.

\item Given $\g,\g' \in C_1$  we must have: $$\be_1  \g\circ\be_1 \g'=\be_1 (\g'\g)$$
$$\g\be_1 \circ \g' \be_1=(\g\g')\be_1 $$
$$\g\be_1 \circ \be_1 \g'= \be_1\g' \circ \g \be_1,$$
whenever these compositions make sense.

\item\label{A} We now define two horizontal compositions of 2-morphisms
$${\begin{pmatrix} & \G'\\\G &
 \end{pmatrix}= \begin{pmatrix} \d^+_2(\G) \be_1 & \G' \\\G  &\be_1 \d_2^-(\G') 
 \end{pmatrix} =\left(\G \be_1 \d_2^-(\G') \right) \be_2\left (\d^+_2(\G) \be_1 \G'\right)} $$and
 $$ {\begin{pmatrix}  \G&\\&\G' 
 \end{pmatrix}=\begin{pmatrix} \G &\be_1 \d_2^+(\G') \\ \d^-_2(\G) \be_1& \G' \end{pmatrix}= \left (\d^-_2(\G) \be_1 \G'\right) \be_2\left(\G \be_1 \d_2^+(\G') \right);}$$ 
and of 3-morphisms:
$${\begin{pmatrix} & J'\\J &
 \end{pmatrix}=
\begin{pmatrix}
 \d^+_2(J)  \be_1  &J'\\J& \be_1  \d_2^-(J') 
\end{pmatrix}=
\left(J \be_1 \d_2^-(J') \right) \be_2\left (\d^+_2(J) \be_1 J'\right)}$$
and 
$${\begin{pmatrix}  J&\\&J' 
 \end{pmatrix}= 
\begin{pmatrix}
 J &\be_1 \d_2^+(J')\\ \d^-_2(J) \be_1 & J'
\end{pmatrix}=
\left (\d^-_2(J) \be_1 J'\right) \be_2\left(J \be_1 \d_2^+(J') \right)}$$ 
It follows from the previous axioms that they are associative. In fact they also define functors $\C_3(x,y) \times \C_3(y,z) \to \C_3(x,z)$, where $\C_3(x,y)$ is the category with objects 2-morphisms $\G$ with $\d^-_2(\G)=x$ and $\d^+_2(\G)=y$ and morphisms the 3-morphisms $J$ with $\d^-_2(J)=x$ and $\d^+_2(J)=y$, and upwards multiplication as composition; {this follows from $7$ and $8$.}

\item {\bf (Interchange 3-cells)} For any two 2-morphisms $\G$ and $\G'$ with $\d_1^+(\G)=\d_1^-(\G')$ a 3-morphism  (called an interchange 3-cell)
$$\begin{pmatrix} & \G'\\\G &
 \end{pmatrix} =\d_3^-(\G \#\G') \ra{ \quad\G\#\G'\quad }  \d_3^+(\G \#\G')= \begin{pmatrix}  \G&\\&\G' 
 \end{pmatrix}$$

\item{\bf (2-functoriality)} For any 3-morphisms $\G_1=\d_3^-(J) \ra{J} \d_3^+(J)=\G_2$ and $\G_1'=\d_3^-(J') \ra{J'} \d_3^+(J')=\G_2'$, with $\d_1^+(J)=\d_1^-(J')$ the following upwards compositions of 3-morphisms coincide:
$$\begin{pmatrix} & \G'_1\\\G_1 &
 \end{pmatrix}  \ra{\quad \G_1\#\G'_1\quad }  \begin{pmatrix}  \G_1&\\&\G_1' 
 \end{pmatrix}\ra{\begin{pmatrix}  J&\\&J '
 \end{pmatrix}} \begin{pmatrix}  \G_2&\\&\G_2' 
 \end{pmatrix}$$
and
$$\begin{pmatrix} & \G'_1\\\G_1 &
 \end{pmatrix}\ra{\begin{pmatrix} & J'\\J & 
 \end{pmatrix}} \begin{pmatrix} & \G'_2\\\G_2 & 
 \end{pmatrix} \ra{\quad \G_2\#\G'_2\quad } \begin{pmatrix}  \G_2&\\&\G_2' 
 \end{pmatrix}.$$
This of course means that the collection $\G \#\G'$, for arbitrary 2-morphisms $\G$ and $\G'$ with $\d^+_1(\G)=\d^-_1(\G')$ defines a natural transformation between the two functors of \ref{A}. Note that by using the interchange condition for the vertical and upwards compositions, we only need to verify this condition for the case when either $J$ or $J'$ is an identity. (This is the way this axiom appears written in \cite{KP,Cr,Be}.)

\item {\bf (1-functoriality)}  For any three 2-morphisms $\g \ra{\G} \f \ra{\G'} \psi$ and $\g ''\ra{\G''}\f''$ with $\d^+_2(\G)=\d^-_2(\G')$ and $\d^+_1(\G)=\d^+_1(\G')= \d^-_1(\G'')$ the following upwards compositions of 3-morphisms coincide:
$$\begin{pmatrix} 
\psi  \be_1 &\G''\\  \G' &\be_1 \g''\\ \G&\be_1 \g''\end{pmatrix}
 \ra{\begin{pmatrix} \G' \# \G''\\ \G \be_1 \g''\end{pmatrix}}
 \begin{pmatrix}  \G'& \be_1 \f''\\ \phi \be_1 & \G''\\  \G&\be_1\g''\end{pmatrix} 
\ra{\begin{pmatrix} \G' \be_1 \f''\\ \G\#\G'' \end{pmatrix}}
\begin{pmatrix}  \G'& \be_1 \f''\\ \G & \be_1 \f''\\ \g\be_1 &\G''
\end{pmatrix} 
 $$
where the 2-morphism components of the 3-morphisms stand for the
corresponding identity 3-morphism, and
$$\begin{pmatrix} 
\psi \be_1 &\G''\\  \G' &\be_1\g''\\ \G&\be_1\g''\end{pmatrix}
 \ra{\begin{pmatrix} \G' \\\G \end{pmatrix} \#\G''}
\begin{pmatrix}  \G'& \be_1 \f''\\ \G & \be_1 \f''\\ \g\be_1 &\G''
\end{pmatrix}. 
 $$
Furthermore an analogous equation holds with left and right whiskering
exchanged.
\end{enumerate}

\begin{Definition}\label{GrayFunctor}
A (strict) Gray functor $F\colon \C \to \C'$ between Gray 3-groupoids $\C$ and $\C'$ is given by maps $\C_i \to \C'_i$ {($i=0,\ldots, 3)$} preserving all compositions, identities, {interchanges} and boundaries, strictly.
\end{Definition}

\subsubsection{From 2-crossed modules to Gray 3-groupoids with a single object}\label{nnn}
Let $\Hc=(L \ra{\de} E\ra{\d} G,\t,\{,\})$ be a 2-crossed module (of groups). We can construct a Gray 3-groupoid $\C$ with a single object out of $\Hc$. We  put $C_0=\{*\}$, $C_1=G$, $C_2=G \times E$ and $C_3=G \times E \times L$. {This construction appears in \cite{KP}, with different conventions, and also in \cite{CCG,BG}, in a slightly different language.}

 As boundaries $\d^\pm_1 \colon C_k \to C_0=\{*\}$, where  $k=1,2,3$, we take the unique possible map. Furthermore:
$$\d^-_2(X,e)=X \an \d^+_2(X,e)=\d(e)^{-1}X. $$
In addition put (as vertical composition): $$(X,e)\be_2 (\d(e)^{-1}X,f)=(X,ef),$$
and also $${X \be_1 (Y,e)=(XY,X \t e) \an (Y,e) \be_1 X=(YX,e)}.$$
Analogously 
$${X \be_1 (Y,e,l)=(XY,X \t e,X \t l) \an (Y,e,l) \be_1 X=(YX,e,l) .}$$

Looking at 3-cells, put 
$${\d^-_3(X,e,l)=(X,e) \an \d^+_3(X,e,l)=(X,\de(l)^{-1}e) } $$
and 
$${\d^-_2(X,e,l)=X \an \d^+_2(X,e,l)=\d(e)^{-1}X.}$$
(Note ${\d^+_2\d^+_3(X,e,l)=\d^+_2(X,\de(l)^{-1}e)=\d(e)^{-1}X=\d^+_2(X,e,l)}$,
since $\d \de=1$.)
As  vertical composition of 3-morphisms we put:
$$\big(X,e,l\big)\be_2\big(\d(e)^{-1}X,f,k\big)=\left (\begin{CD}\big(\d(e)^{-1}X,f,k\big)\\ \big(X,e,l\big) \end{CD}\right)= \big(X,ef,(e \t ' k)l\big),$$
and as upwards composition  of 3-morphisms  we  put:
$$\big(X,e,l\big)\be_3\big(X,\de(l)^{-1} e,k\big)=\big(X,e,lk\big)$$
The   vertical and upwards compositions of 2-cells define a 2-groupoid since $\big(\de \colon G\times  L \to  G\times E, \t'\big)$ is a crossed module of groupoids \cite{BHS}; see the comments after definition \ref{2cmlg}. Recall that {the category of crossed modules of groupoids and the category of 2-groupoids are equivalent}.

Let us now define the interchange 3-cells. We can see that:
$$\begin{pmatrix} & (Y,f)\\(X,e) & \end{pmatrix}= \left(XY,e \big(\d(e)^{-1} X \big) \t f \right)$$
and 
$$\begin{pmatrix}(X,e) &\\ & (Y,f) \end{pmatrix}
 =\left(XY,\big(X \t f\big) e\right).$$
We therefore take:
\begin{equation}(X,e)\#(Y,f)=\Big(XY,e \big(\d(e)^{-1} X \big) \t f,   e \t' \left\{e^{-1},X \t f\right\}^{-1}\Big). \end{equation}
Note 
\begin{multline*}
\de\left(e \t' \left\{e^{-1},X \t f\right\}^{-1}\right)^{-1}e\big(\d(e)^{-1} X \big) \t f \\=e e^{-1} (X \t f) e \big(\d(e)^{-1}X \big)\t f^{-1} e^{-1} e\big(\d(e)^{-1} X \big) \t f=(X \t f)e.
\end{multline*}

It is easy to see that:
\begin{align*}
\begin{pmatrix} & (Y,f,l)\\(X,e,k) & \end{pmatrix}&=\begin{pmatrix} (\d(e)^{-1}X Y, \d(e)^{-1}X  \t f, \d(e)^{-1}X \t l \\ (XY,e,k)\end{pmatrix}\\&=\Big(XY,e \big(\d(e)^{-1} X \big) \t f, \big(e \t' \d(e)^{-1} X \t l\big)k\Big)
\end{align*}
and 
\begin{align*}
\begin{pmatrix}(X,e,k) &\\ & (Y,f,l) \end{pmatrix}
 &=\begin{pmatrix}
(X\d(f)^{-1}Y,e,k)\\(XY,X \t f,X \t l)
  \end{pmatrix}\\ &=\Big(XY,(X\t f)e,((X \t f) \t' k)X \t l\Big).
\end{align*}
To prove condition 13. of the definition of a Gray 3-groupoid (2-functoriality) we must prove that (for each $X \in G$, $e,f \in E$ and $k,l \in L$):
\begin{multline*}
e \t' \{e^{-1}, X \t f\}^{-1}((X \t f) \t' k)X \t l\\= \big(e \t' \d(e)^{-1} X \t l\big)k (\de(k)^{-1} e) \t' \{e^{-1} \de(k),X \t (\de(l)^{-1} f)\}^{-1}.
\end{multline*}
or, by using the fact that $(\de\colon L \to E, \t')$ is a crossed module:\begin{multline}\label{referr}
e \t' \{e^{-1}, X \t f\}^{-1}((X \t f) \t' k)X \t l\\= \big(e \t' \d(e)^{-1} X \t l\big)  e \t' \{e^{-1} \de(k),X \t (\de(l)^{-1} f)\}^{-1} k.
\end{multline}

For $l=1$ this is equivalent to:
$$
e \t' \{e^{-1}, X \t f\}^{-1}((X \t f) \t' k)=e \t' \{e^{-1} \de(k),X \t  f\}^{-1} k.
$$
or:
$$((X \t f)\t' k^{-1}) e \t'\{e^{-1},X \t f\}=k^{-1} e \t'\{e^{-1} \de(k),X\t f\}$$
which follows from equation (\ref{imppp}) and the definition of $e \t' l=l\{\de(l^{-1}),e\}$. {Note that $\d \circ \delta=1_L$.}
For $k=1$ equation (\ref{referr}) is the same as:
\begin{equation*}
e \t' \{e^{-1}, X \t f\}^{-1}X \t l= \big(e \t' \d(e)^{-1} X \t l\big)  e \t' \{e^{-1} ,X \t (\de(l)^{-1} f)\}^{-1},
\end{equation*}
or
\begin{equation*}
(X \t l^{-1} )e \t' \{e^{-1}, X \t f\}=  e \t' \{e^{-1} ,X \t (\de(l)^{-1} f)\}\big(e \t' \d(e)^{-1} X \t l^{-1}\big).
\end{equation*}
This can be proved as follows, by using equation (\ref{Porter})
\begin{align*}
 &e \t' \{e^{-1} ,X \t (\de(l)^{-1} f)\}\big(e \t' \d(e)^{-1} X \t l^{-1}\big)
\\&=\left( \de(X \t l^{-1})e \t' \{e^{-1},X \t f  \}\right) \left( e \t '\{e^{-1},X \t \de(l)^{-1}\} \right) \big(e \t' \d(e)^{-1} X \t l^{-1}\big)
\\&=(X \t l^{-1}) e \t' \{e^{-1},X \t f  \} (X \t l)\left( e \t '\{e^{-1},X \t \de(l)^{-1}\} \right) \big(e \t' \d(e)^{-1} X \t l^{-1}\big)
\end{align*}
where we have used the fact that $(\de\colon L \to E,\t')$ is a crossed module. Now note, by using $6.$ of the definition of a 2-crossed module:
\begin{align*}
 &\left( e \t '\{e^{-1},X \t \de(l)^{-1}\} \right) \big(e \t' \d(e)^{-1} X \t l^{-1}\big)\\&=e \t '\left ( \{X \t \de(l^{-1}), e^{-1}\}^{-1} (X \t l^{-1})    \right)
\\&=(e \de(X \t l^{-1}))\t '\left ( \{X \t \de(l), e^{-1}\} \right) e \t' \left (X \t l^{-1}    \right), \quad \textrm{by equation } (\ref{invv})
\\&=(e \t' X \t l^{-1} ) e\t' \left ( \{X \t \de(l), e^{-1}\} \right), {\textrm{ since }  (\de\colon L \to E,\t')  \textrm{ is a crossed module}}\\
&=e \t'\big ( (X \t l^{-1} ) \{X \t \de(l), e^{-1}\} \big)
\\&=X \t l^{-1}, \textrm{ by definition of } e \t'l=l\{\de(l^{-1}),e\}.
\end{align*}
The general case of equation (\ref{referr}) follows from $k=1$ and $l=1$ cases by the interchange law for the upwards and vertical compositions.

Let us now prove 1-functoriality (condition 14. of the definition of a Gray 3-groupoid). The first condition is equivalent to:
$$(ef) \t' \left\{f^{-1}e^{-1},X \t g \right\}^{-1}=(ef)\t' \left\{ f^{-1},\d(e)^{-1}X \t g\right \}^{-1}e\t'\left\{e^{-1}, X \t g\right\}^{-1} .$$
This follows directly from equation (\ref{imppp}). The second condition is equivalent to:
$$e\t'\{e^{-1}, X \t f X \t g\}^{-1}=e \t' \{e^{-1}, X\t f\}^{-1} \big((X \t f) e\big)\t' \{e^{-1},X \t g\}^{-1} $$
which follows from equation (\ref{Porter}).

{We have therefore proved that any 2-crossed module $\H$ defines a Gray 3-groupoid $\C(\H)$, with a single object. This process is reversible: a  Gray 3-groupoid  $C$ together with an object $x \in C_0$ of it defines a  2-crossed module; see \cite{KP,Be,BG}.}

\subsubsection{{Example: (finite type) chain complexes}}\label{l3cc}

{Suppose ${A=\{A_n,\d_n=\d\}}_{n \in \mathbb{Z}}$ is a chain complex of finite dimensional vector spaces, such that the set of all $n$ for which $A_n$ is not the trivial vector space is finite. (Chain complexes like this will be called of finite type.)}  This construction is analogous to the one in \cite{KP}.
 
{
Let us then construct a Lie 2-crossed module $$\GL(A)=\left(\GL^3(A) \ra{\a}  
\GL^2(A) \ra{\b} \GL^1(A),\t,\{,\}\right)$$ out of $A$. If $A$ is not of finite type, then the same construction will still yield a  crossed module, albeit of infinite dimensional Lie groups.}

  The  group $\GL^1(A)$ is given by all invertible chain maps $f\colon A \to A$, with composition as product. This is a Lie group, with Lie algebra $\gl^1(A)$ given by all chain maps $A \to A$, with bracket given by the usual commutator of chain maps. We also denote the algebra of all chain maps $A\to A$ with composition as product as $\hom^1(A)$.

Recall that a homotopy is given by a degree-1 map $s\colon A \to A$.  Let $\hom^2(A)$ denote the vector space of 1-homotopies. Likewise, we  define an $n$-homotopy as being a {degree-$n$ map $b\colon A \to A$, and denote the vector space of $n$-homotopies as $h^{n+1}(A)$.} Notice also that we have a complex $\{h^n(A), \d_n'\}$, where $\d'_n(b)=\d b-(-1)^{n}b \d$, for each $b \in \hom^n(A)$. {Note $h^1(A)=\hom^2(A)$ and $\hom^1(A)\subset h^0(A)$.}

We define the $*$ product of two 1-homotopies as:
$$s*t=s+t+s\d t+st\d$$
This defines an associative product in $\hom^2(A)$. Even though $\hom^2(A)$ is not an algebra, considering commutators this yields a Lie algebra $\gl^2(A)$, with commutator:
$$[s,t]=st\d+s\d t-ts\d-t\d s .$$
Note that $\gl^2(A)$  is the Lie algebra of the Lie group $\GL^2(A)$ of invertible elements of $\hom^2(A)$, the identity of this latter group being the null homotopy.

It is easy to see  that the map $\b\colon \hom^2(A) \to \hom^1(A)$ such that
$$\b(s)=1+\d s +s\d $$ {respects the products.} This thus defines a Lie group morphism $\b\colon \GL^2(A) \to \GL^1(A)$, the differential form of which is given by the Lie algebra map $\b'\colon \gl^2(A) \to \gl^1(A)$, where $$\b'(s)=\d s +s \d.$$

There is also a left action of $\GL^1(A)$ on $\GL^2(A)$ by automorphisms given  by 
$$ f \t s=fsf^{-1}. $$ Its differential form is given by the left action of $\gl^1(A)$ on $\gl^2(A)$ by derivations such that:
$$f \t s=fs-sf. $$
It is easy to see that we have defined a pre-crossed module of Lie groups and of Lie algebras. Moreover, this yields a crossed module if the chain complex is of length 2, as we will see below.

By definition a 3-track will be an element of $\hom^3(A)=h^2(A)/\d'(h^3(A))$, therefore it will be a 2-homotopy up to a 3-homotopy. Considering the sum of 3-tracks, defines an abelian Lie group $\GL^3(A)$ whose Lie algebra ${\gl}^3(A)$ is 
given by the vector space of 3-tracks, with trivial commutator.

The {(well defined)}  map $\a\colon \GL^3(A) \to \GL^2(A)$ such that $$\a(b)=-\d b+ b\d$$ is a group morphism.
 We also have a left  action of $\GL^1(A)$ on $\GL^3(A)$ by automorphisms defined as $$f\t a=faf^{-1}. $$ This defines a complex of Lie groups acted on by $\GL^1(A)$:
$${\GL}^3(A) \ra{\a} \GL^2(A) \ra{\b} \GL^1(A). $$ Its differential form is:
$${\gl}^3(A) \ra{\a'} \gl^2(A) \ra{\b'} \gl^1(A), $$
where $\gl^3(A)$  is the vector space $\hom^3(A)$ with trivial commutator, and $\a'=\a$. Note that $\gl^1(A)$ acts on the left on $\gl^3(A)$ by derivations as $f \t b=fb-bf$.

To define a 2-crossed module we now need to specify the Peiffer lifting. We can see that given $s,t \in \GL^2(A)$ we have
\begin{align*}
s*t*s^{-1}&=\b(s)t\b(s)^{-1}-\d st\b(s)^{-1}+st\d+st\d s^{-1} \d \\
        &=\b(s)t\b(s)^{-1} +\a(st)\b(s)^{-1}\\
       &=(\b(s)t\b(s)^{-1})*(\a(st)\b(s)^{-1})\\
       &=(\b(s)t\b(s)^{-1})*\a(st\b(s)^{-1})
\end{align*}
In particular the Peiffer pairing is:
$$\left <s,t\right>=\left(\b(s)\t t\right) *\a\left(st\b(s)^{-1}\right)* \left(\b(s)\t t^{-1}\right) .$$
This can still be simplified.
Let $a \in {\GL}^3(A)$. We can see that $$t*\a(a)*t^{-1}=\a\left(a\b(t)^{-1} \right)=\a\left(a\right)\b(t)^{-1} .$$
This follows by applying the penultimate equation, noting that:
\begin{align*}
t* \a(a)*t^{-1}&=(\b(t)\a(a)\b(t)^{-1})*\a\big(t\a(a)\b(t)^{-1}\big) \\
&=\a(\b(t)a\b(t)^{-1})*\a\big(t\a(a)\b(t)^{-1}\big) \\
&=\a(\b(t)a +t\a(a)\big)\b(t)^{-1}\\
&=\a(\d t a+t a\d  + a  )\b(t)^{-1}\\
&=\a(a)\b(t)^{-1}
\end{align*}
The Peiffer pairing thus simplifies to:
$$\left<s,t\right> = \a\left( st \b(t)^{-1}\b(s)^{-1}\right),$$
and we thus have the following candidate for the role of the Peiffer lifting:
$$ \{s,t\}=st\b(t)^{-1}\b(s)^{-1},$$
where $s,t\in \GL^2(A)$. Its differential form is:
$$ \{s,t\}=st,$$
where $s,t\in \gl^2(A)$.

Routine calculations prove that we have indeed defined 2-crossed modules of Lie groups and of Lie algebras. The fact we are considering 3-tracks (2-homotopies up to 3-homotopies), instead of simply 2-homotopies, is used several times to prove this.

\subsubsection{Example: the automorphism 2-crossed module of a crossed module}

Let $\Gc=(\d \colon E \to G, \t)$ be a Lie group crossed module. Let us build the differential 2-crossed module associated with the automorphism 2-crossed module of $\Gc$. In the case of crossed modules of groups,  the construction of this 2-crossed module appears in \cite{BG,RS,N}, in the latter in the language of crossed squares. The extension to crossed modules of Lie groups is straightforward.

Let $\lG=(\d\colon \le \to\lg,\t)$ be the differential crossed  module associated to $\Gc$. Let us then construct a 2-crossed module of Lie algebras
$$\left (\gl^3(\lG) \to \gl^2(\lG) \to \gl^1(\lG),\t,\{\,\}\right).$$

The Lie algebra $\gl^1(\lG)$ is given by all chain maps $f=(f_2,f_1)\colon \lG \to \lG$, which, {termwise,} are Lie algebra derivations: 
$$f_2([u,v])=[f_2 (u),v]+[u,f_2(v)], \fo u,v \in \le$$
and
$$f_1([x,y])=[f_1 (x),y]+[x,f_1(y)], \fo x,y \in \lg,$$
satisfying additionally:
$$f_2(x \t v)=f_1(x) \t v+x \t f_2(v), \fo x \in \lg \an v \in \le. $$
The Lie algebra structure is given by the termwise commutator of derivations.

The Lie algebra $\gl^2(\lG)$ is given by all pairs $(x,s)$, where $s\colon \lg \to \le $ is a linear map such that:
$$s([x,y])=x \t s(y)-y \t s(x),$$ 
(in other words $s\colon \lg\to \le$ is a derivation, and we put 
$s \in \Der(\lg,\le)$)  
 and $x \in \lg$. {The Lie algebra structure on $\gl^2(\lG)=\lg \ltimes \Der(\lg,\le)$ is given by a semidirect product, as we now explain.}

The commutator of two derivations $s,t\in  \Der(\lg,\le)$ is:
$$[s,t]=s\d t-t \d s. $$
It is easy to see that this is also a derivation. The crossed module relations are used several times to prove this. (This would not be true if a pre-crossed module was used.)

There exists a left action of $\gl^1(\lG)$ on the Lie algebra of derivations $s\colon \lg \to\le$ given by:
$$(f_1,f_2) \t s=f_2 s-sf_1. $$
We also have a Lie algebra map $q\colon \Der(\lg,\le) \to \gl^1(\lG)$ given by $$q(s)=\d s+s\d.$$ This defines a differential crossed module. 

There is another crossed module of Lie algebras that can be constructed from $\lG$. This is provided by the map $q'=(q'_1,q'_2)\colon \lg \to \gl^1(\lG)$ which associates to each $x\in \lg$ the inner derivation $f_x\colon \lg \to \lg$ such that $f_x(y)=[x,y],$ for each $y \in \lg$, and the derivation $\le \to \le$ such that $v\mapsto x \t v$. The action of $\gl^1(\lG)$ on $\lg$ is $(f_1,f_2)\t x=f_1(x)$, where  $x \in \lg$.

{In particular, we also have an action of $\lg$ on $\Der(\lg,\le)$, provided by the map $q'\colon \lg \to \gl^1(\lG)$ and the already given action of $\gl^1(\lG)$ on $\Der(\lg,\le)$.} {Therefore we can put a  Lie algebra structure on  $\gl^2(\lG)$  given by the semidirect product $\lg\ltimes \Der(\lg,\le)$.} In particular 
$$[(x,s),(y,t)]=([x,y],x \t t-y \t s+[s,t]).$$ 
The boundary map $\b'\colon \gl^2(\lG) \to \gl^1(\lG)$ is:
$${\b'(a,s)=q'(a)+q(s).}$$ 
By the above, this is a Lie algebra map. We also define $f\t (a,s)=(f\t a, f \t s)$, which defines a differential  pre-crossed module. The Peiffer pairing is given by:
$$ \left<(x,s),(y,t)\right>=-\left(\d s(y),F_{s(y)}\right)$$
Here, given $e \in\le$, the map $F_e\colon \lg \to \le $ is $F_e(x)=x \t e$.

The Lie algebra $\gl^3(\lG)$ is given by $\le$. The boundary map $\a'\colon \gl^3(\lG) \to \gl^2(\lG)$ is  $\a'(e)=(\d e,F_e)$.

 The Peiffer lifting is defined as: $$\{(x,s),(y,t)\}=-s(y).$$
Therefore: 
$$\{\a'(e),\a'(f)\}=\{( \d e,F_e),( \d f,F_f)\} =-F_e(\d f)=-\d f \t e=-[f,e]=[e,f].$$
Furthermore the action of $\gl^1(\lG)$ on $\le$ is $(f_1,f_2) \t e=f_2(e)$.

The (rest of the) straightforward proof that this defines a differential 2-crossed module is left to the reader. {Note that given a Lie group $G$ we can define a differential crossed module $(\id \colon \lg \to \lg, \ad)$. The  automorphism 2-crossed module  of it is exactly the differential 2-crossed module associated to the 2-crossed module of Example \ref{trivialexample}.  }

\section{The thin fundamental Gray 3-groupoid of a smooth manifold}
Let $M$ be a smooth manifold. {We denote $D^n=[0,1]^n$.}

\begin{Definition}[$n$-path]
Let $n$ be a positive integer. An $n$-path is given by a smooth map $\a\colon D^n=D^1\times D^{n-1}\to M$ for which there exists an $\e>0$ such that $\a(x_1,x_2,\ldots x_n)=\a(0,x_2,\ldots x_n)$ if $x_1 \leq \e$, and analogously for any other face of $D^n$, of any dimension. We will abbreviate this condition as saying that $\a$ has a product structure close to the boundary of the $n$-cube. We also suppose that $\a(0\times D^{n-1})$ and $\a(1\times D^{n-1})$ each consist of just a single point.
\end{Definition}

Given an $n$-path and an $i\in \{1,\ldots, n\}$ we can define $(n-1)$-paths $\d^-_i(\a)$ and $\d^+_i(\a)$ by restricting $f$ to $D^{i-1} \times \{0\} \times D^{n-i}$ and $D^{i-1} \times \{1\} \times D^{n-i}$. Note that $\d^\pm_1(\a)$ are necessarily constant $(n-1)$-paths.
Given two $n$-paths $\a$ and $\b$ with $\d^+_i(\a)=\d^-_i(\b)$ we consider the obvious concatenation $\a\be_i b$, which given the product structure condition of $\a$ and $\b$ is also an $n$-path; see examples below.

\subsection{1-Tracks (the rank-1 homotopy relation)} \label{1-Tracks}
Note that a 1-path is given by a smooth path $\gamma\colon [0,1] \to M$ such that there exists an $\epsilon >0$ such that $\g$ is constant in $[0,\e] \cup [1-\e,1]$, which can be abbreviated by saying that each end point of $\g$ has a sitting instant (we are using the terminology of \cite{CP}). Given a 1-path $\g$, define the source and target of $\g$ as $\d^-_1(\g)=\g(0)$ and $\d^+_1(\g)=\g(1)$, respectively. 

Given two 1-paths $\g$ and $\f$ with $\d^+_1(\g)=\d^-_1(\f)$, their concatenation $\g\f=\g\be_1\f$ is  the usual one
$$(\g\f)(t)=\left \{ \begin{CD} \g(2t), \textrm{ if }t \in [0,1/2] \\ \f(2t-1), \textrm{ if } t \in [1/2,1]\end{CD} \right.$$ 
The fact that any 1-path has sitting instants at its end points implies that the concatenation of two 1-paths is again a 1-path.

Similarly, a 2-path $\G$ is given by a smooth map $\G\colon [0,1]^2 \to M$ such that  there exists an $\epsilon >0$ for which:
\begin{enumerate}
\item $\G(t,s)=\G(0,0)$ if $0 \leq t \leq \epsilon$ and $s \in [0,1]$,
\item $\G(t,s)=\G(1,0)$ if  $1-\epsilon \leq t \leq 1$ and $s \in [0,1]$,
\item $\G(t,s)=\G(t,0)$ if $0 \leq s \leq \epsilon$ and $t \in [0,1]$
\item $\G(t,s)=\G(t,1)$ if  $1-\epsilon \leq s \leq 1$ and $t \in [0,1]$.
\end{enumerate}

The following definition appeared in \cite{CP}. See also \cite{FMP1,FMP2,MP,M}.
\begin{Definition}[Rank-1 homotopy]
Two 1-paths $\f$ and $\g$ are said to be rank-1 homotopic  (and we write $\f\cong_1 \g$) if there exists a 2-path $\G$ such that:
\begin{enumerate}
\item $\d^-_2(\G)=\g$ and $\d^+_2(\G)=\f$.
\item {$\rank(\Dc \G(v)) \leq 1, \forall v \in [0,1]^2.$}
\end{enumerate}
Here $\Dc$ denotes derivative.
\end{Definition}
Note that if $\g$ and $\f$ are rank-1 homotopic, then they have the same initial and end-points.  Given the product structure condition on 2-paths, it follows that rank-1 homotopy is an equivalence relation. Given a 1-path $\g$, the equivalence class to which it belongs is denoted by $[\g]$, or simply $\g$, when there is no ambiguity.

We denote the set of $1$-paths of $M$ by $S_1(M)$. The quotient of $S_1(M)$ by the relation of thin homotopy is denoted by $\S_1(M)$. We call the elements of $\S_1(M)$  1-tracks. 

It is easy to prove that the concatenations of $1$-tracks (defined in the obvious way from the concatenation of 1-paths) together with the source and target maps $\sigma, \tau \colon \S_1(M) \to M$, defines a groupoid $\Sc_1(M)$ whose set of morphisms is  $\S_1(M)$ and whose set of objects is $M$. (For details see \cite{CP}.)

\begin{Definition}
Let $* \in M$ be a base point. The group $\pi_1^1(M,*)$ is defined as being the set of  1-tracks $[\g]\in \S_1(M)$ starting and ending at $*$, with the group operation being the concatenation of {1-tracks}. 
\end{Definition}

\subsection{Strong and laminated 2-tracks}
\subsubsection{Strong 2-Tracks (the strong rank-2 homotopy relation)}\label{2-Tracks}
\begin{Definition}[Strong rank-2 homotopy]\label{thin2}
Two 2-paths $\G$ and $\G'$ are said to be strong rank-2 homotopic  (and we write $\G\cong_2^s \G'$) if there exists a 3-path $J\colon D^3 \to M$ such that:
\begin{enumerate}

\item We have $\d_3^-(J)=\G$ and $\d_3^+(J)=\G'$.
\item The restrictions $\d^\pm_2(J)$ restrict to rank-1 homotopies $\d^\pm_2(\G) \to \d^\pm_2(\G')$.
\item  {$\rank ( \Dc {J}(v)) \leq 2$ for any $v \in {[0,1]}^3$.}
\end{enumerate}
\end{Definition}
Due to the fact that any $3$-path has a product structure close to its boundary, it follows that strong rank 2-homotopy is an equivalence relation. 

We denote by $S_2(M)$ the set of all 2-paths of $M$. The quotient of $S_2(M)$ by the relation of strong rank-2 homotopy is denoted by $\S_2^s(M)$. We call the elements of $\S_2^s(M)$ strong 2-tracks.

{This notion of strong rank-2 homotopy was used in \cite{MP,SW2,SW3,BS,FMP1,FMP2,M}, and it behaves very nicely with respect to 2-dimensional holonomy based on a crossed module. For this reason, it is too strong for our purposes in this article. Therefore we define now a weaker version of rank-2 homotopy.}
\subsubsection{Laminated 2-Tracks (the laminated rank-2 homotopy equivalence relation)} 

\begin{Definition}[Laminated rank-2 homotopy]\label{laminated}
Two 2-paths $\G$ and $\G'$ are said to be  laminated rank-2 homotopic (and we write $\G\cong_2^l \G'$) if there exists a 3-path $J\colon D^3 \to M$ (say $J(t,s,x))$ such that:
\begin{enumerate}
\item We have $\d_3^-(J)=\G$ and $\d_3^+(J)=\G'$, in other words $J(t,s,0)=\G(t,s)$ and $J(t,s,1)=\G'(t,s)$ for each $s,t\in[0,1]$.
\item The restrictions $\d^\pm_2(J)$ (in other words $J(t,0,x)$ and $J(t,1,x)$) restrict to rank-1 homotopies $\d^\pm_2(\G) \to \d^\pm_2(\G').$
\item {$\rank ( \Dc {J}(v)) \leq 2$ for any $v \in {[0,1]}^3$.}
\item For each $0<s,x <1$, {at least}  one of the following conditions holds {(up to a set of Lebesgue measure zero)}:
\begin{enumerate}
\item {\bf Laminatedness} For either {$\zeta =s$}  or {$\zeta =x$} the following condition holds for all $t\in [0,1]$:  $$\rank(\Dc_{(t,\zeta)} J(t,s,x)) \leq 1.$$ 
\item {\bf Path space thinness:} There exist non-zero constants $a$ and $b$ such that 
$$a\frac{\d}{\d s} J(t,s,x)+b\frac{\d}{\d x} J(t,s,x)=0$$
for each $t\in [0,1]$.
\end{enumerate}
\end{enumerate}
\end{Definition}
Once again, since any $3$-path has a product structure close to {its boundary, it follows that} laminated rank-2 homotopy is an equivalence relation.

The quotient of $S_2(M)$, the set of 2-paths of $M$, by the relation of laminated rank-2 homotopy is denoted by $\S_2^l(M)$. We call the elements of $\S_2^l(M)$ laminated 2-tracks.
{Note that the boundaries $\d^\pm_i\colon S_2(M) \to S_{i-1}(M)$ {descend to boundaries $\d^\pm_i\colon \S_2^l(M) \to \S_{i-1}(M)$ and $\d^\pm_i\colon \S_2^s(M) \to \S_{i-1}(M)$; here $i=1,2$.}}
\subsubsection{Identity 2-tracks}\label{identities}
There exists a  map $S_1(M) \to S_2(M)$ sending a path $\g$ to the 2-path $\id(\g)$ (frequently written simply as $\g$) such that  $\id(\g)(t,s)=\g(t),$ for each $s,t \in [0,1].$ It descends to maps $\id\colon \S_1(M) \to \S_2^s(M)$ and $\id\colon \S_1(M) \to \S_2^l(M)$.

\subsubsection{A technical lemma}
The following lemma will be needed for proving the {consistency} of the vertical composition of laminated and strong 2-tracks. It generalises {Lemma} 52 of \cite{FMP1}.

\begin{Lemma}\label{coher}
Let $f\colon \d (D^3) \to M$ be a smooth map  such that:
\begin{enumerate}
 \item The restriction of $f$ to any face of $\d D^3$ has a product structure close to the boundary of it.
\item The restrictions $f(0,s,x)$ and $f(1,s,x)$ are constant.
\item We have {$\rank (\Dc f(v))\leq 1, \forall v \in \d D^3$.}
\end{enumerate}
Then $f$ extends to a map $g\colon D^3 \to M$ defining a laminated rank-2 homotopy connecting $\d^-_3(f)=\d^-_3(g)$ and  $\d^+_3(f)=\d^+_3(g)$.
\end{Lemma}
\begin{Proof}
Let $D^3=\{(t,s,x):-1\leq t,s,x \leq 1\}$ and $S^2$ be its boundary; a smooth manifold with corners. According to the proof of {Lemma} 52 of \cite{FMP1}, the map $f\colon S^2 \to M$  factors as $f=p \circ \f$, where $\f\colon S^2 \to N$ is a smooth map, {$N$} being a contractible manifold, and $p\colon N \to M$ {being a smooth local diffeomorphism. In particular conditions 1,2 and 3 hold for $\f$.}

{Choose a contraction $c\colon N \times [0,1] \to N$  of $N$ to a point $*$ of it. We can suppose that, at each end of $[0,1]$, $c$ has a sitting instant, in other words, that there exists a positive $\e$ such that $c(x,t)=x$ if $t \in [0,\e]$ and $c(x,t)=*$ if $t \in [1-\e,1]$, where $x \in N$.}

Consider the map $g\colon D^3 \to M$ defined as:
$${g(t,s,x)=p\left( c\left(\f\Big(\frac{(t,s,x)}{|(t,s,x)|}\Big), 1-|(t,s,x)| \right)\right) }$$
where $|(t,s,x)|=\max\left(|t|,|s|,|x|\right)$. {The map $g$} is smooth and has a product structure close to the boundary of $D^3$. All this follows from the product structure condition for $\phi$ on the faces of $S^2$ and the sitting instant condition on $c$, which take care of the arguments where $(t,s,x) \mapsto |(t,s,x)|$ is not smooth.

The map $g\colon D^3 \to M$ extends $f\colon \d D^3 \to M$, and therefore we now only need to prove that it is a laminated rank-2 homotopy.  Note that it follows trivially that $g$ is a rank-2 homotopy, as in the proof of of {Lemma} 52 of \cite{FMP1}.

Define $|(s,x)|=\max(|s|,|x|)$. Let $(s,x)\in [-1,1]^2$. If $(s,x)$ is a point where $|(s,x)|$ is smooth then either $\dx |{(s,x)}|=0$ or $\ds |{(s,x)}|=0$. Suppose $\dx |{(s,x)}|=0$, from which it follows that $\dx |{(t,s,x)}|=0$, for each $t \in [0,1]$. Let us see that
\begin{equation}\label{rrr}\rank(\Dc_{(t,x)} g(t,s,x))\leq 1, 
\fo  t\in [0,1].\end{equation}  Given a $t \in [0,1]$, we either  have $\dt |(t,s,x)|=0$ or not. In the first case (\ref{rrr}) follows from the fact  {$\rank(\Dc f(v))\leq 1, \forall v \in S^2$} together with  $\dx |{(t,s,x)}|=0$, for each $t \in [0,1]$. In the second case (\ref{rrr}) follows from the fact that $f$ is constant when $t=1$ or $t={-1}$, which makes $g$ depend only on $t$ in a neighbourhood of $(t,s,x)$.

 The same argument is valid when $\ds |{(s,x)}|=0$. 
The remaining points $(s,x)$ have measure zero.
\end{Proof}

\subsubsection{Vertical composition of 2-tracks}\label{vcomp2}
Recall that we can vertically compose any two 2-paths $\G$ and $\G'$ with $\d^+_2(\G)=\d^-_2(\G')$. Denote it by $\G\be_2\G'$, and represent it graphically as: $$\G\be_2\G'=\begin{CD} \d_2^+(\G') \\ @AA\G'A\\\d_2^-(\G')=\d_2^+(\G) \\@AA \G A\\ \d_2^-(\G)  \end{CD}=\begin{CD} \G'\\\G \end{CD} .$$

Suppose $[\G]$ and $[\G']$ are  (laminated or strong 2-tracks) such that $\d^+_2([\G])=\d^-_2([\G])$. Choose a rank-1 homotopy $H$ connecting $\d^+_2(\G)$ and $\d^-_2(\G')$. Define $[\G]\be_2 [\G']=[\G\be_2 H \be_2 \G']$.
By using {Lemma} \ref{coher} {we obtain} the following, not entirely trivial, result:
\begin{Lemma}
The vertical composition of (laminated or strong) 2-tracks is well defined (does not depend on any of the choices made).
\end{Lemma}
See \cite{BH1,BHS,HKK,FMP1,BHKP} for similar constructions.

It is easy to see that the vertical  composition of strong and laminated 2-tracks is associative. In the laminated case, we will need to use the path-space thinness condition in definition \ref{laminated}. We have:
\begin{Proposition}
 The vertical composition of strong or laminated 2-tracks  defines  categories with morphisms $\S_2^s(M)$ and $\S_2^l(M)$, respectively, and objects $\S_1(M)$. The source and target maps are $\d^-_2$ and $\d^+_2$. The identities are as in \ref{identities}. 
\end{Proposition}

\subsubsection{Whiskering  2-tracks by 1-tracks}\label{1whis}
Let $\G$ be a 2-path. Let also $\g$ be a 1-path, such that $\d^+_1(\G)=\d^-_1(\g).$ The right whiskering of $\G$ with $\g$ is by definition the 2-path $\G \be_1 \g\doteq \G \be_1 \id(\g)$; see \ref{identities}. We analogously define left whiskering $\g' \be_1 \G$ if $\d^+_1(\g)=\d^-_1(\G)$, and whiskering of $n$-paths by 1-paths for arbitrary $n$.

{It is easy to show that {these} whiskerings descend to an action of the groupoid $\Sc_1(M)$ on the categories $\S_2^l(M)$ and $\S_2^s(M)$. The {main} part of the proof is to show that, in the laminated case,
$[\G] \be_1 [\g]\doteq [\G \be_1 \g]$ does not depend on the representatives $\G$ and $\g$ chosen. Suppose we have a rank 1-homotopy $H$ connecting $\g_1$ and $\g_2$ and a laminated rank-2 homotopy $J$ connecting $\G_1$ and $\G_2$. Then $(J \be_1 \g_1) \be_3 ( \id(\G_2)\be_1 \id(H))$ is a laminated rank-2 homotopy connecting $\G_1 \be_1 \g_1 $ and $\G_2 \be_1 \g_2$. Here $\id(\G_2)(t,s,x)=\G_2(t,s)$ and $\id(H)(t,s,x)=H(t,x)$, where $t,s,x \in [0,1]$.}

\subsubsection{Horizontal compositions of strong 2-tracks}\label{ftg}

Let $\G$ and $\G'$ be two 2-paths with ${\d^+_1(\G)=\d^-_1(\G')}$. We can consider the obvious horizontal concatenation $\G\be_1\G'$, denoted by $ \G\G'$.
There are however two other natural ways to define the horizontal composition of $\G$ and $\G'$. These are:
$$ {\G\be_1^-\G'}=\big(\G \be_1 \d^-_2(\G')\big)\be_2 \big(\d^+_2(\G)\be_1 \G'\big)=\left (\begin{CD} \quad &\G'\\  \G & \quad \end{CD}\right)$$
and
$$ {\G\be_1^+\G'}=\big(\d^-_2(\G)\be_1 \G'\big) \be_2 \big(\G \be_1 \d^+_2(\G')\big) =\left(\begin{CD}  &\G \quad\\ &\quad \G'  \end{CD}\right)$$
It is easy to see that  these three horizontal compositions descend to the quotient $\S_2^s(M)$ of $S_2(M)$ under strong rank-2 homotopy and they all coincide. In addition the interchange law between the horizontal and vertical compositions holds. In fact  we have:
\begin{Theorem}
The horizontal and vertical composition of strong 2-tracks {defines} a 2-groupoid $\Sc_2^s(M)$ with objects given by the set of points of $M$, 1-morphisms given by $\S_1(M)$ and 2-morphisms by $\S_2^s(M)$.
\end{Theorem}
See \cite{FMP1,FMP2,MP,SW2} for details. The definition of a 2-groupoid appears, for example,  in \cite{HKK}. For related constructions see \cite{BH1,BHS,BH3,BHS,HKK}.
\subsubsection{Horizontal compositions of laminated 2-tracks and the interchange 3-track}\label{i3c}

Let us now look at the behaviour of the horizontal composition under the relation of laminated rank-2 homotopy. We can see that {$\be_1$ does not descend to the quotient. However $\be_1^-$ and $\be_1^+$ do descend, even though they do not coincide}.

Given 2-paths $\G$ and $\G'$ with $\d^+_1(\G)=\d^-_1(\G')$, let $J=\G \# \G'$ be the 3-path whose typical slices 
 as {$x$ varies} appear in figure \ref{ash}. We have put $\g_1=\d^+_2(\G)$, $\g_0=\d^-_2(\G)$, $\g'_x(t)=\G'(t,x)$, $\G_x'(t,s)=\G'(t,xs)$ and {$\hat{\G_x'}(t,s)=\G'(t,x+s(1-x))$}, where $t,s,x \in [0,1]$. Note $\G'_x \be_2 \hat{\G'_x}=\G'$, for each $x \in [0,1]$.
\begin{figure}
\centerline{\relabelbox 
\epsfysize 2.3cm
\epsfbox{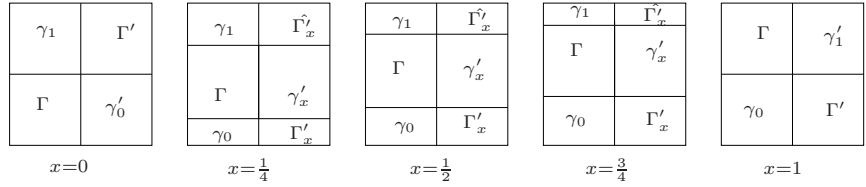}
\relabel{a}{$\scriptstyle{\G}$}
\relabel{b}{$\scriptstyle{\G'}$}
\relabel{d}{$\scriptstyle{\g_1}$}
\relabel{c}{$\scriptstyle{\g'_0}$}
\relabel{a1}{$\scriptstyle{\G}$}
\relabel{b1}{$\scriptstyle{\g'_x }$}
\relabel{d1}{$\scriptstyle{ \g_1}$}
\relabel{e1}{$\scriptstyle{\g_0 }$}
\relabel{f1}{$\scriptstyle{\G_x' }$}
\relabel{g1}{$\scriptstyle{\hat{\G_x'}}$}
\relabel{a2}{$\scriptstyle{\G}$}
\relabel{b2}{$\scriptstyle{\g'_x }$}
\relabel{d2}{$\scriptstyle{ \g_1}$}
\relabel{e2}{$\scriptstyle{\g_0 }$}
\relabel{f2}{$\scriptstyle{\G_x' }$}
\relabel{g2}{$\scriptstyle{\hat{\G_x'}}$}
\relabel{a3}{$\scriptstyle{\G}$}
\relabel{b3}{$\scriptstyle{\g'_x }$}
\relabel{d3}{$\scriptstyle{ \g_1}$}
\relabel{e3}{$\scriptstyle{\g_0 }$}
\relabel{f3}{$\scriptstyle{\G_x' }$}
\relabel{g3}{$\scriptstyle{\hat{\G_x'}}$}
\relabel{a4}{$\scriptstyle{\G}$}
\relabel{b4}{$\scriptstyle{\G'}$}
\relabel{e4}{$\scriptstyle{\g_0 }$}
\relabel{g4}{$\scriptstyle{\g_1'}$}
\relabel{A}{$\scriptstyle{x=\frac{1}{4}}$}
\relabel{B}{$\scriptstyle{x=\frac{1}{2}}$}
\relabel{C}{$\scriptstyle{x=\frac{3}{4}}$}
\relabel{D}{$\scriptstyle{x=0}$}
\relabel{E}{$\scriptstyle{x=1}$}
\endrelabelbox}
\caption{\label{ash} Slices of $\G\#\G'$ for  $x=0$, $x=1/4$, $x=1/2$, $x=3/4$ and $x=1$.}
\end{figure}

This good (see below) 3-path $\G \# \G'$ is well defined up to reparametrisations in the $x$-direction. Below (see subsection \ref{3t}) we define an equivalence relation on 3-paths which solves this ambiguity.

Note that $\d^-_3(\G\# \G')=\G \be_1^- \G'$ and $\d^+_3(\G\# \G')=\G \be_1^+ \G'.$ This will give us the interchange 3-cells; see \ref{dg3g}.

\subsection{3-Tracks}\label{3t}
\begin{Definition}[Good 3-path]\label{good} {A 3-path $(t,s,x) \mapsto J(t,s,x)$ is called good if $\d^\pm_2(J)$ each are independent of $x$.}
\end{Definition}

\begin{Definition}[rank-3 homotopy (with laminated boundary)]\label{3strong}
We  will \\say that two good 3-paths $J$ and $J'$ are  rank-3 homotopic (with laminated boundary), and we write $J\cong_3 J'$ if there exists a 4-path  $(t,s,x,u)\in D^4 \mapsto W(t,s,x,u)\\\in M$ such that:
\begin{enumerate}
\item We have $\d_4^-(W)=J$ and $\d_4^+(W)=J'$, in other words $W(t,s,x,0)=J(t,s,x)$ and $W(t,s,x,1)=J'(t,s,x)$, where $t,s,x,u \in [0,1]$.
\item The restriction $W(t,1,x,u)$ is independent of $x$ and defines a  rank-1 homotopy connecting $W(t,1,x,0)$ and $W(t,1,x,1)$ (each independent of $x$, therefore identified with paths in $M$), and the same for the restriction $W(t,0,x,u)$.
\item The restriction $W(t,s,0,u)$ defines a laminated rank-2 homotopy connecting $J(t,s,0)$ and $J'(t,s,0)$, and analogously for $W(t,s,1,u)$.
\item {For each $v\in [0,1]^4$ we have $\rank\left (\Dc W(v) \right)\leq 3$.}
\end{enumerate}
\end{Definition}
We denote  $\S_3(M)$ as being {the set} of all good 3-paths up to  rank-3 homotopy (with laminated boundary). The elements of $\S_3(M)$ will be called 3-tracks. 

Notice that the {interchange 3-track } $\G\#\G'$ of two 2-tracks with $\d^+_1(\G)=\d^-_1(\G')$ gives rise to a well defined element of $\S_3(M)$. In addition, all boundaries $\d^\pm_i,i=1,2,3$ of good 3-paths descend to maps $\d^\pm_1\colon \S_3(M) \to M$, $\d^\pm_2\colon \S_3(M) \to \S_1(M)$ and $\d^\pm_3\colon \S_3(M) \to \S_2(M)$, which are the maps needed to get the structure of a Gray 3-groupoid (\ref{dg3g}.)

The following lemma will be useful later. The proof is {achieved} by using a filling argument very similar to the proof of {Theorem} A of \cite{BH2}.
\begin{Lemma}\label{Refercoher1}
Two good 3-paths $J$ and $J'$ are rank-3 homotopic (with laminated boundary) if, and only if, there exists a 4-path $W\colon D^4\to M$ satisfying the conditions 1,3,4 of  definition \ref{3strong} but with condition $2$ replaced by: 
\end{Lemma}
\begin{enumerate}
\item [2'] {The restriction $W(t,1,x,u)$  defines a laminated rank-2 homotopy which connects the 2-paths $W(t,1,x,0)$ and $W(t,1,x,1)$}, and the analogous condition holds for the restriction $W(t,0,x,u)$.
\end{enumerate}
\begin{Proof}
One of the implications $(2 \implies 2')$ is immediate. Let us prove the reciprocal. We will just discuss how to deal with the condition on the $s=1$ face, since the other face is dealt with analogously.

If $W$ is a homotopy connecting $J$ and $J'$ as in the statement of the lemma, we substitute $W$ by  $W\be_2 V$, where $V$ is defined in the following way: consider a smooth retraction {$r\colon D^3 \to \d D^3\setminus \{x=1\}$}. Define a smooth function $U\colon \left  ( \d D^3\setminus \{x=1\}\right)\times [0,1] \to M$, as:
\begin{enumerate}
 \item $U(t,0,x,u)=W(t,1,x,u)$, the right hand side is a laminated rank-2 homotopy. 
\item In {${\big(\d(D^3) \setminus (\{x=1\} \cup \{s=0\})\big)\times [0,1]}$ we put $U(t,s,x,u)=W(t,1,0,u)$. Note $W(t,1,0,u)$ is a rank-1 homotopy.}
\end{enumerate}
Finally, let {$V(t,s,x,u)=U(r(t,s,x),u))$}. Then $W\be_2 V$ is a rank-3 homotopy with laminated boundary, and will give us, after a very minor adjustment {(making use of the path-space thinness condition)}, a rank-3 homotopy with laminated boundary connecting $J$ and $J'$.
\end{Proof}

\subsubsection{Vertical composition of 3-tracks}
Let $J$ and $J'$ be good 3-paths with $\d^+_2(J)=\d^-_2(J')$. Recall that we can perform their vertical composition $J \be_2 J'$. If $[J]$ and $[J']$ are such that $\d^+_2([J])=\d^-_2([J'])$, then choosing a rank-1 homotopy $H$ connecting $\d^+_2(J)$ and $\d^-_2(J')$ permits us to put $[J]\be_2 [J']\doteq [J \be_2 \id(H) \be_2 J']$, where $\id(H)(t,s,x)=H(t,s)$, for each $t,s,x \in [0,1]$. By using the same argument as in \ref{vcomp2}, we can see that this vertical composition of 3-tracks is well defined.

\subsubsection{Upwards composition of 3-tracks}

Let $J$ and $J'$ be good 3-paths with $\d_3^+(J)=\d_3^-(J')$. Consider the composition $J\be_3 J'$, called upwards composition. 
Suppose that we have $\d_3^+(J)\cong_2^l\d_3^-(J')$. Choose a laminated rank-2 homotopy $H$ connecting $\d_3^+(J)$ and $\d_3^-(J')$, and put: 
$[J]\be_3 [J']\doteq [J\be_3 H \be_3 J']$.
The proof of the following essential lemma (very inspired by \cite{BH2,BHS}) will be very similar to  the proof of {Lemma} \ref{Refercoher1} and will make use of it.
\begin{Lemma}
This upwards composition is well defined in $\S_3(M)$.
\end{Lemma}\begin{Proof}
 Suppose $J\cong_3 J_1$ and $J'\cong_3 J_1'$, and choose  rank-3 homotopies with laminated boundary $A$ and $A'$ yielding {these} equivalences. Choose also  laminated rank-2 homotopies with  $H$ connecting $\d_3^+(J)$ and $\d_3^-(J')$ and $H_1$ connecting $\d_3^+(J_1)$ and $\d_3^-(J'_2)$.

Consider the map $W\colon D^4 \to M$ defined in the following way: we fill the  $x=1$ and $x=0$ faces of $D^4$ with the laminated rank-2 homotopies $\d_3^-(A')$ and $\d_3^+(A)$. Then we fill the $u=0$ and $u=1$ faces of $D^4$ with $H$ and $H'$. The boundary of the $s=0$ face of $D^4$ will define a rank-1 homotopy $U(t,v)$ connecting $H(t,0,0)$ with itself. Explicitly {it is given by} the following concatenation of rank-1 homotopies:
\begin{multline*}H(t,0,0) \ra{H(t,0,x)} H(t,0,1)\ra{A'(t,0,0,u)} A'(t,0,0,1) =H'(t,0,1) \\{\ra{H'(t,0,1-x)} H'(t,0,0) \ra{A(t,0,1,1-u)}H(t,0,0)}
\end{multline*}
By applying  {Lemma} \ref{coher}, the $s=0$ face of $D^4$ can therefore be filled with a laminated rank-2 homotopy. To extend $W$ to the rest of $D^4$ put $W(t,s,x,u)=W(t,r(s,x,u))$, where $r$ is a smooth retraction of $D^3$ onto $\d D^3 \setminus \{s=1\}$. 

Then $A \be_3 W \be_3 A'$ {satisfies} the conditions of {Lemma} \ref{Refercoher1}, which implies there exists  {a rank-3 homotopy with laminated boundary} connecting $J\be_3 H \be_3 J'$  and $J_1\be_3 H_1 \be_3 J_1'$.
\end{Proof}
\subsubsection{Whiskering  3-tracks by 1-tracks}
The treatment is entirely similar to what was presented in \ref{1whis}.

\subsubsection{The fundamental thin Gray 3-groupoid $\Sc_3(M)$ of a smooth manifold $M$}\label{thingray}
{Combining all of the above we have:}
\begin{Theorem}
 {Let $M$ be a smooth manifold. The sets of 1-tracks, laminated 2-tracks and 3-tracks can be arranged into a Gray 3-groupoid $$\Sc_3(M)=\left(M,\S_1(M),\S_2^l(M),\S_3(M)\right)$$ whose set of objects is $M$.}
\end{Theorem}
\begin{Proof}
We need to verify the conditions in \ref{dg3g}. Conditions 1. to 3. are trivial. Note that $\d_2^\pm$, applied to good 3-paths, can be naturally regarded as mapping to 1-paths.

The difficult bit (existence of compositions) of conditions 4. to 10. are already proved, and all the rest follows {straightforwardly} by the definition of laminated and rank-3 homotopy (with laminated boundary), as in the construction in \cite{FMP1,FMP2,MP,SW2}. We have already proved the existence of an interchange 3-cell; \ref{i3c}. The fact that it verifies conditions 13. and 14. follows since both sides of each equation can be connected by rank-3 homotopies with a laminated boundary.
\end{Proof}

\section{Three-dimensional holonomy based on a Gray 3-groupoid}
Fix a smooth manifold $M$. We will make use of Chen's definition of differential forms in the smooth space of smooth paths in  $M$, as well as iterated integrals of differential forms; see \cite{Ch}. For conventions see the Appendix. 

The main result of this section is:
\begin{Theorem}\label{Main}
Consider  a 2-crossed module  {$\Hc=(L\ra{\de} E \ra{\d} G, \t, \{,\})$ with associated differential  2-crossed module $(\ll\ra{\de} \le \ra{\d} \lg, \t, \{,\})$.} Let $M$ be a smooth manifold. Consider differential forms  $\w \in \A^1(M,\lg)$, $m \in \A^2(M,\le)$ and $\theta \in \A^3(M,\ll)$  such that $\de(\theta)=\M$ and $\d(m)=\W$, where  $\W=d\w+[\w,\w]=d\w+\frac{1}{2}\w \wedge^\ad \w$ and $\M=dm+\w \wedge^\t m$ denote the curvature of $\w$ and 2-curvature 3-form of the pair $(m,\w)$; see \cite{FMP1,FMP2,BS,SW2,SW3}. 

Then we can define a (smooth) strict Gray 3-groupoid functor (definition \ref{GrayFunctor}) $$\stackrel{(\w,m,\theta)}{\H}\colon \Sc_3(M) \to \C(\Hc)  ,$$
where $\C(\Hc)$ is the  Gray 3-groupoid constructed from $\Hc$ (and whose sets of objects and 1,2 and 3-morphisms are smooth manifolds.)
\end{Theorem}
{The explicit description  of $\H$ appears in \ref{US}. At the level of 1-paths $\H$ coincides with the usual holonomy of a non-abelian 1-form, whereas for crossed modules it yields the 2-dimensional holonomy of \cite{SW2,BS,FMP2}.}
\begin{Definition}
 Given $\w,m$ and $\theta$ as above, the 3-curvature 4-form $\Theta$ of $(\w,m,\theta)$ is given by
$$\Theta=d\theta + \w \wedge^\t \theta-m \wedge^{\{,\}} m ;$$
see the appendix for this notation. {Here $m \wedge^{\{,\}} m$ is the antisymmetrisation of $6\{m, m\}$.}
\end{Definition}

\subsection{{One-dimensional} holonomy based on a Lie group}\label{1dim}
Let $G$ be a Lie group with Lie algebra $\lg$ and $A\colon \R \to \lg$ be a smooth map. Let $F^A(t_0,t)$ be the solution of the differential equation:
$$\dt F^A (t_0,t)=F^A(t_0,t)A(t), \textrm{ with } F^A(t_0,t_0)=1_G$$
It thus follows that $$\frac{\d}{\d t_0} F^A (t_0,t)=-A(t_0)F^A(t_0,t).$$
Moreover
 $$ F^A(t_0,t+t')=F^A(t_0,t)F^A(t,t').$$
More generally, suppose that $(t,s) \in \R^2 \mapsto A_s(t)\in \lg$ is a smooth map. It is well known (and not difficult to prove) that:
$$\ds F^{A_s}(t_0,t)=\int_{t_0}^t  F^{A_s}(t_0,t') \ds A_s(t')F^{A_s}(t',t) dt'.$$

Let $\w$ be a 1-form in the manifold $M$. Let $\g\colon [0,1] \to M$ be a smooth curve in $M$, and $A$ be given by {$\g^*(\w)=A (\dot \g)$}. Put $g^\w_\g(t_0,t)=F^A(t_0,t)$. Note $g^\w_\g(t_0,t+t')=g^\w_\g(t_0,t)g^\w_\g(t,t')$ and {$g^\w_\g(t,s)=g^\w_\g(s,t)^{-1}$.}

Suppose $s \in I \mapsto \g_s$ is a smooth one-parameter family of smooth curves in $M$. In other words the map  $\G\colon [0,1]^2 \to M$ such that $\G(t,s)=\g_s(t)$ for each $s,t \in [0,1]$ is smooth. We have:
\begin{align*}
&\frac{\d}{\d s} g^\w_{\g_s}(a,b)\\&=\int_{a}^b g^\w_{\g_s}(a,t)\frac{\d}{\d s} \w\left(\frac{\d}{\d t} \g_s(t)\right)g^\w_{\g_s}(t,b)dt\\
&=\int_{a}^b g^\w_{\g_s}(a,t)
 \left(
d\w\left(\frac{\d}{\d s} \g_s(t), \frac{\d}{\d t} \g_s(t)\right)
+\frac{\d}{\d t}\w\left(\frac{\d}{\d s} \g_s(t)\right)\right)
g^\w_{\g_s}(t,b)dt\\
&=g^\w_{\g_s}(a,b)\int_{a}^b g^\w_{\g_s}(b,t)
 \left(
d\w\left(\frac{\d}{\d s} \g_s(t), \frac{\d}{\d t} \g_s(t)\right)
+\frac{\d}{\d t}\w\left(\frac{\d}{\d s} \g_s(t)\right)\right)
g^\w_{\g_s}(t,b)dt
\end{align*} 
Now note that (integrating by parts):
\begin{align*}
 &\int_{a}^b (g^\w_{\g_s}(t,b))^{-1}
 \frac{\d}{\d t}\w\left(\frac{\d}{\d s} \g_s(t)\right)g^\w_{\g_s}(t,b)dt
\\&\quad=\left.(g^\w_{\g_s}(t,b))^{-1}\w\left(\frac{\d}{\d s} \g(t,s)\right)(g^\w_{\g_s}(t,b))\right|_{t=a}^{t=b}\\
&\quad\quad+\int_{a}^b (g^\w_{\g_s}(t,b))^{-1} \left[\w\left(\frac{\d}{\d s} \g_s(t)\right),\w\left(\frac{\d}{\d t} \g_s(t)\right)
\right]
 g^\w_{\g_s}(t,b)dt
\end{align*}
And therefore, putting $\W=d\w+[\w,\w]$ as being the curvature of $\w$ (and where $[\w,\w](X,Y)=[\w(X),\w(Y)]$), we {get} the following {very well known} lemma:
\begin{Lemma}\label{OneHol}
Let $\w\in \A^1(M,\lg)$ be a $\lg$-valued 1-form in $M$. We have:
\begin{multline}
\frac{\d}{\d s} g^\w_{\g_s}(a,b)=-g^\w_{\g_s}(a,b)\int_{a}^b (g^\w_{\g_s}(t,b))^{-1} \W\left(\frac{\d}{\d t} \g_s(t),\frac{\d}{\d s} \g_s(t)\right)
 g^\w_{\g_s}(t,b)dt\\+g^\w_{\g_s}(a,b)\left (\left.(g^\w_{\g_s}(t,b))^{-1}\w\left(\frac{\d}{\d s} \g(t,s)\right)(g^\w_{\g_s}(t,b))\right|_{t=a}^{t=b}\right).
\end{multline}
{This can also be written as:}
\begin{multline}\label{zxc}
\frac{\d}{\d s} g^\w_{\g_s}(a,b)^{-1}=g^\w_{\g_s}(a,b)^{-1}\int_{a}^b g^\w_{\g_s}(a,t) \W\left(\frac{\d}{\d t} \g_s(t),\frac{\d}{\d s} \g_s(t)\right)
 g^\w_{\g_s}(t,a)^{-1}dt\\-g^\w_{\g_s}(a,b)^{-1}\left (\left. g^\w_{\g_s}(a,t)\w\left(\frac{\d}{\d s} \g(t,s)\right)(g^\w_{\g_s}(a,t))^{-1}\right|_{t=a}^{t=b}\right).
\end{multline}
\end{Lemma}

{We thus arrive at the following well known result; see \cite{CP,MP,FMP1,BS}.}
\begin{Corollary}\label{1thin}
One dimensional holonomy based on a Lie group is invariant under rank-1 homotopy. More precisely, if $\g$ and $\g'$ are rank-1 homotopic 1-paths then $g^\w_\g(0,1)=g^\w_{\g'}(0,1)$.
\end{Corollary}

\subsubsection{A useful lemma of Baez and Schreiber}\label{Useful}
Suppose that the group $G$ has a left action $\t$, {by automorphisms, on} the Lie group $E$, with Lie algebra $\le$.  
Let $*,*'\in M$. We define $P(M,*,*')$ to be the space of all smooth paths $\g\colon [0,1]\to M$ that start {at $*$ and finish at} $*'$.

Let ${\U}$ be an open set of $M$.
Let {$f\colon (t,x) \in [0,1] \times \U \mapsto \g_x(t)=f_t(x)\in M$} define a smooth map $F\colon \U \to P(M,*,*')$, which is the same as saying that $f$ is smooth. Note that we have  {$\g_x(0)=*$ and $\g_x(1)=*',\forall x \in \U$}. Let $A$ be a $(n+1)$-form in $M$ with values in $\le$. Let $\w$ be a 1-form in $M$ with values in $\lg$. Put $D_\w A=d A+\w\wedge^\t A$, the exterior covariant derivative of $A$.

Consider  the twisted iterated integral:
\begin{equation}
\oint_0^b g^\w_{\g_x} \t f^*( A)=\int_0^b g^\w_{\g_x}(0,t) \t \left(\i_{\dt }f^*( A)\right)dt \in \A^n(\U,\le).
\end{equation}
By using the {first equation} of the previous lemma we have
\begin{align*}d( g^\w_{\g_x} (0,t))&= {- \left(\int_{0}^t (g^\w_{\g_x}(0,t'))\t^{\ad} \i_{\dt}f^*(\W)dt'\right) g^\w_{\g_x} (0,t)+ g^\w_{\g_x}(0,t)f^*(\w).}
\end{align*}
Here the map {$(t,x) \mapsto g^\w_{\g_x}(0,t) \in G$} should be seen (through the action $\t$) as taking values in the vector space of linear maps $\le \to \le$.
Applying equation (\ref{Ch1})  {from the Appendix} it thus follows:
\begin{multline*}
d\oint_0^b g^\w_{\g_x} \t f^*( A) \\=-\oint_0^b g^\w_{\g_x} \t D_\w f^*(A)-\oint_0^b \left(g^\w_{\g_x} \t^\ad f^*(\W)\right)  *^\t\left( g^\w_{\g_x} \t  f^*(A)\right)+f_b^*\left(g^\w_{\g_x} \t A \right)
\end{multline*}
We have used the identity $\i_X(\a\wedge \b)=\i_X(\a)\wedge \b+J(\a)\wedge \i_X(\b)$, valid for any two forms $\a$ and $\b$. Recall {$J(\a)=(-1)^n \a$}, where $n$ is the degree of $\a$.

We  define the following form in the loop space $P(M,*,*')$:
$$\oint_\w A=\oint \left(g_\g^{\w}(0,t)\right) \t A,$${and analogously for iterated integrals.
 We are following the notation 
(but not the conventions) of \cite{BS}.}
In other words: $${F^*\left (\oint_\w A\right )=\oint_0^1 g^\w_{F(x)}(0,t) \t f^*( A)dt ,}$$
for each plot $F\colon \U \to P(M,*,*')$, the map $f\colon I \times \U \to M$ being $f(t,x)=F(x)(t)$, where $x \in \U$ and $t\in I$.

{We thus have the following very useful lemma which appeared in \cite{BS}.}
\begin{Lemma}[Baez-Schreiber]\label{BSL}
{Let the Lie group $G$ act on the Lie group $E$ by automorphisms. Consider a smooth manifold $M$. Let $A$ be a $(n+1)$-form in $M$ with values in $\le$. Let $\w$ be a 1-form in $M$ with values in $\lg$. We have:}
\begin{equation}
d \oint_\w A=-\oint_\w D_\w A -\oint_\w \W*^\t A,
\end{equation}
where $$D_\w A=d A+\w \wedge^\t A$$
is the exterior covariant derivative of $A$ with respect to $\w$. 
\end{Lemma}

\subsection{Two-dimensional holonomy based on a pre-crossed module}
Suppose that $\Gc=(\d\colon E \to G, \t)$ is  a (as usual Lie) pre-crossed module. Let $\lG=(\d\colon \le \to \lg, \t)$ be the associated differential pre-crossed module. Let $\w \in \A^1(M,\lg)$ be a $\lg$-valued smooth {1-form in $M$.} Let $m\in \A^2(M,\le)$ be an arbitrary $\le$-valued 2-form in $M$. (Later we will put the restriction $\d(m)=\W$, where $\W=d\w+[\w,\w]=d \w+\frac{1}{2}\w\wedge^\t w$ is the curvature of $\w$).  {See the Appendix for notation.}

Given a smooth map $\G\colon (t,s) \in [0,1]^2 \mapsto 
 \g_s(t) \in M$, defining therefore a smooth map {(plot) $s\in [0,1] \stackrel{\hat{\G}}{\mapsto} \g_s \in P(M,*,*')$,} define $e_\G^{(\w,m)}(s_0,s)$  as the solution of the differential equation:
\begin{equation}\label{zxcv}
\frac{d}{d s} e_\G^{(\w,m)}(s_0,s)= e_\G^{(\w,m)}(s_0,s) \int_{0}^{1} g^\w_{\g_s}(0,t) \t m\left(\frac{\d}{\d t} \g_s(t),\frac{\d}{\d s} \g_s(t)\right)
 dt,
\end{equation}
with initial condition $e_\G^{(\w,m)}(s_0,s_0)=1_E$.
In other words
$$\frac{d}{ds}  e_\G^{(\w,m)}(s_0,s)=e_\G^{(\w,m)}(s_0,s){\hat{\G}}^*\left(\oint_\w m\right);$$
see \ref{Useful} {and \ref{forms}} for this notation.
Note that {(by using {Lemma} \ref{BSL})} the curvature of the 1-form  $\oint_\w m$ (in the path space $P(M,*,*')$) is 
$$d \oint_\w m+\left[\oint_\w m,\oint_\w m\right]=-\oint_\w \M-\oint_\w \W *^\t m +\left[\oint_\w m,\oint_\w m\right],$$
where $\M=dm+\w\wedge^\t m$ is defined as being the 2-curvature $3$-form of $(\w,m)$.
By using {Lemma} \ref{OneHol}, it thus follows that:
\begin{Theorem}
 Let $M$ be a smooth manifold. Suppose  $\Gc=(\d\colon E \to G, \t)$ is  a Lie pre-crossed module, with associated differential {pre}-crossed module $\lG=(\d\colon \le \to \lg,\t)$. Let $\w \in \A^1(M,\lg)$ and $m\in \A^2(M,\le)$. Let $\M=dm+\w \wedge^\t m$ be the 2-curvature of ${(\w,m)}$.
Suppose that $(t,s,x) \in {[0,1]^3} \mapsto \G_x(t,s)={\G^s(t,x)}=\g_s^x(t)=J(t,s,x)\in M$ is smooth. {Suppose   that $\g_s^x(0)=*$ and $\g_s^x(1)=*'$ for each $s,x \in [0,1]$.} {Let $\hat{J}\colon [0,1]^2 \to P(M,*,*')$ be the associated plot, {$(s,x) \mapsto \g_s^x $},}
  We have:
\begin{multline*}
\frac{d}{dx}e_{\G_x}^{(\w,m)}(0,a)\\=\int_0^a e_{\G_x}^{(\w,m)}(0,s){\hat{J}}^*\left(\oint_\w \M+\oint_\w \W *^\t m-\left [\oint_\w m, \oint_\w m \right]\right)\left(\ds,\dx\right)e_{\G_x}^{(\w,m)}(s,a)ds\\ + e_{\G_x}^{(\w,m)}(0,a){\hat{\G^a}}^*
\left (\oint_\w m\right)\left(\dx\right)- {\hat{\G^0}}^*\left(\oint_\w m\right) \left(\dx\right)e_{\G_x}^{(\w,m)}(0,a).
\end{multline*}
Explicitly:
\begin{align*}
&\frac{d}{dx}e_{\G_x}^{(\w,m)}(0,a)\\&=\quad\int_0^a\int_0^1 e_{\G_x}^{(\w,m)}(0,s)  
\left (
g^\w_{\g_s^x}(0,t)\t \M\left(\frac{\d}{\d t} \g_s^x(t),\frac{\d}{\d s} \g_s^x(t),\frac{\d}{\d x} \g_s^x(t)\right)\right)e_{\G_x}^{(\w,m)}(s,a)dtds\\ 
&\quad+
\int_0^a\int_0^1  e_{\G_x}^{(\w,m)}(0,s)\left(\int_0^t g^\w_{\g_s^x}(0,t') {\t^{\ad}}\W\left(\frac{\d}{\d t'} \g_s^x(t'),\frac{\d}{\d s} \g_s^x(t')\right)dt'\right)g^\w_{\g_s^x}(0,t)\t\\ & \quad\quad\quad\quad\quad\quad\quad m\left(\frac{\d}{\d t} \g_s^x(t),\frac{\d}{\d x} \g_s^x(t)\right)e_{\G_x}^{(\w,m)}(s,a)dtds
\\
&\quad-
\int_0^a\int_0^1  e_{\G_x}^{(\w,m)}(0,s)\left(\int_0^t g^\w_{\g_s^x}(0,t') {\t^{\ad}}\W\left(\frac{\d}{\d t'} \g_s^x(t'),\frac{\d}{\d x} \g_s^x(t')\right)dt'\right)g^\w_{\g_s^x}(0,t)\t\\&\quad\quad\quad\quad\quad m\left(\frac{\d}{\d t} \g_s^x(t),\frac{\d}{\d s} \g_s^x(t)\right)e_{\G_x}^{(\w,m)}(s,a)dtds
\\&\quad-
\int_0^a e_{\G_x}^{(\w,m)}(0,s) \Big[ \int_0^1 g^\w_{\g_s^x}(0,t')\t  m\left(\frac{\d}{\d t'} \g_s^x(t'),\frac{\d}{\d s} \g_s^x(t')\right)dt' ,\\
& \quad\quad\quad\quad\quad \int_0^1 g^\w_{\g_s^x}(0,t)\t m\left(\frac{\d}{\d t} \g_s^x(t),\frac{\d}{\d x} \g_s^x(t)\right)dt\Big]e_{\G_x}^{(\w,m)}(s,a)ds\\
&\quad+e_{\G_x}^{(\w,m)}(0,a)
\int_0^1 g^\w_{\g_a^x}(0,t)\t m\left(\frac{\d}{\d t} \g_a^x(t),\frac{\d}{\d x} \g_a^x(t)\right)dt\\
&\quad-\left(\int_0^1 g^\w_{\g_0^x}(0,t)\t m\left(\frac{\d}{\d t} \g_0^x(t),\frac{\d}{\d x} \g_0^x(t)\right)dt\right)e_{\G_x}^{(\w,m)}(0,a)
\end{align*}
\end{Theorem}
{Now note that for any smooth functions $f$ and $g$ we have $[\int_0^1 f(t)dt,\int_0^1 g(t')dt']=\int_0^1 [f(t),\int_0^t g(t')dt']dt+\int_0^1 [\int_0^tf(t')dt',g(t)]dt$,} and use it in the fourth term. Consider also the definition of  the Peiffer commutators $\left <u,v\right >=[u,v]-\d(u)\t v, \wh u,v \in \le $. We obtain the following result:
\begin{Corollary}\label{HOL2}
 {Under the conditions of the previous theorem, and, furthermore, assuming that} $\d(m)=\W$, we have:
\begin{multline*}
\dx e_{\G_x}^{(\w,m)}(0,a)\\=\int_0^a e_{\G_x}^{(\w,m)}(0,s){\hat{J}}^*\left(\oint_\w \M-\oint_\w m *^{\left <,\right >} m\right)\left(\ds,\dx\right)e_{\G_x}^{(\w,m)}(s,a)ds\\+e_{\G_x}^{(\w,m)}(0,a){\hat{\G^a}}^*
\left (\oint_\w m\right)\left(\dx\right)- {\hat{\G^0}}^*\left(\oint_\w m\right) \left(\dx\right)e_{\G_x}^{(\w,m)}(0,a),
\end{multline*}
{which can also be written as:}
\begin{multline}\label{comp2}
\dx e_{\G_x}^{(\w,m)}(0,a)^{-1}\\=-e_{\G_x}^{(\w,m)}(0,a)^ {-1}\int_0^a e_{\G_x}^{(\w,m)}(0,s){\hat{J}}^*\left(\oint_\w \M-\oint_\w m *^{\left <,\right >} m\right)\left(\ds,\dx\right)e_{\G_x}^{(\w,m)}(s,0)ds
\\-{\hat{\G^a}}^*\left (\oint_\w m\right)\left(\dx\right)e_{\G_x}^{(\w,m)}(0,a)^{-1} +e_{\G_x}^{(\w,m)}(0,a)^{-1}{\hat{\G^0}}^*\left(\oint_\w m\right) \left(\dx\right).
\end{multline}
Explicitly (looking at the first expression):
\begin{align*}
&\frac{d}{dx}e_{\G_x}^{(\w,m)}(0,a)\\&=\quad\int_0^a\int_0^1 e_{\G_x}^{(\w,m)}(0,s)\left( g^\w_{\g_s^x}(0,t)\t \M\left(\frac{\d}{\d t} \g_s^x(t),\frac{\d}{\d s} \g_s^x(t),\frac{\d}{\d x} \g_s^x(t)\right)\right)e_{\G_x}^{(\w,m)}(s,a)dtds\\ 
&\quad-
\int_0^a\int_0^1  e_{\G_x}^{(\w,m)}(0,s)\Big<\int_0^t g^\w_{\g_s^x}(0,t') \t {m}\left(\frac{\d}{\d t'} ,\frac{\d}{\d s} \g_s^x(t')\right)dt',\\ & \quad\quad\quad\quad\quad\quad\quad g^\w_{\g_s^x}(0,t)\t m\left(\frac{\d}{\d t} \g_s^x(t),\frac{\d}{\d x} \g_s^x(t)\right)e_{\G_x}^{(\w,m)}(s,a)dt\Big>ds
\\
&\quad+
\int_0^a\int_0^1  e_{\G_x}^{(\w,m)}(0,s)\Big<\int_0^t g^\w_{\g_s^x}(0,t') {\t m}\left(\frac{\d}{\d t'} \g_s^x(t'),\frac{\d}{\d x} \g_s^x(t')\right)dt',\\ & \quad\quad\quad\quad\quad\quad\quad g^\w_{\g_s^x}(0,t)\t m\left(\frac{\d}{\d t} \g_s^x(t),\frac{\d}{\d s} \g_s^x(t)\right)e_{\G_x}^{(\w,m)}(s,a)dt\Big>ds
\\
&\quad+e_{\G_x}^{(\w,m)}(0,a)
\int_0^1 g^\w_{\g_a^x}(0,t)\t m\left(\frac{\d}{\d t} \g_a^x(t),\frac{\d}{\d x} \g_a^x(t)\right)dt\\
&\quad-\left(\int_0^1 g^\w_{\g_0^x}(0,t)\t m\left(\frac{\d}{\d t} \g_0^x(t),\frac{\d}{\d x} \g_0^x(t) \right)dt\right)e_{\G_x}^{(\w,m)}(0,a).
\end{align*}
{Note that if $(\d\colon E \to G,\t)$ is a crossed module then all terms involving Peiffer commutators vanish. {This will be of prime importance now and later.}}
\end{Corollary}
Using this last  result it follows that:
\begin{Corollary}\label{Thin2}
The two dimensional holonomy based on a pre-crossed module $(\d\colon E \to G,\t)$ is invariant under laminated rank-2 homotopy. More precisely, if $\G,\G'$ are 2-paths {which} are laminated rank-2 homotopic, and $m\in \A^2(M,\le)$ and $\w\in \A^1(M,\lg)$ are such that $\d(m)=\W=d \w+[\w,\w]$ then {$e_\G^{(\w,m)}(0,1)=e_{\G'}^{(\w,m)}(0,1)$}. Two dimensional holonomy based on a crossed module (where, by definition, the Peiffer commutators vanish) is invariant under {strong} rank-2 homotopy.
\end{Corollary}
\begin{Proof}
{If $J(t,s,x)=\g^x_s(t)$ is a laminated rank-2 homotopy then the right hand side of the {last} equation in the previous lemma vanishes. If $J$ is a strong rank-2 homotopy, then in principle only the first term vanishes, but this is compensated by the fact that the Peiffer pairing is zero.}
\end{Proof}

\subsubsection{{The {Non-abelian} Green Theorem}}\label{greensection}
Fix a pre-crossed module $(\d\colon E \to G, \t)$, with associated differential 2-crossed module $(\d\colon \le \to \lg, \t)$.
{The {next result} follows directly from the unicity theorem for ordinary differential equations, together with equations (\ref{zxc}) and (\ref{zxcv}).}
\begin{Lemma}[Non-abelian Green Theorem]\label{Green}
Let $\G$ be a 2-path in $M$. Let  $m\in \A^2(M,\le)$ and $\w\in \A^1(M,\lg)$ be such that $\d(m)=\W=d \w+[\w,\w]$. We then have:
$$ \d\left(e^{(\w,m)}_\G(s,s')\right)^{-1}g^\w_{\g_s}=g^\w_{\g_{s'}}, \fo s,s'\in [0,1].$$ 
{Recall} that we put $\G(t,s)=\g_s(t),\fo s,t \in [0,1]$ and $g^\w_{\g_s}\doteq g^\w_{\g_s}(0,1)$.
\end{Lemma}

{Combining corollaries \ref{1thin} and \ref{Thin2} with the {Non-abelian} Green  {Theorem} it follows that {(recall the notation of \ref{thingray}):}}
\begin{Theorem}{Given $(\w,m)$ as in the previous lemma, 
the assignments: $$\g \in \S_1(M) \mapsto  g_\g^\w \in G, \an \G\in \S_2^l(M) \mapsto \left(g_{\d_1^-(\G)}^\w,e_\G^{(\w,m)}\right) \in G \times E, $$
where $ g_\g^\w\doteq g_\g^\w(0,1)$ and $e_\G^{(\w,m)}\doteq e_\G^{(\w,m)}(0,1)$ satisfy all the axioms for a morphisms of Gray 3-groupoids $\Sc_3(M) \to \C(\H)$ which do not involve 3-morphisms (is short it defines a morphism of sesquigroupoids, \cite{S}).}
\end{Theorem}
\begin{Proof}
What is left to prove is entirely analogous to the proof of {Theorem} 39 of \cite{FMP1}.
\end{Proof}
\begin{Remark}
 {If $(\d\colon E \to G,\t)$ is a crossed module then the assignments}:
 $$\g \in \S_1(M) \mapsto  g_\g^\w \in G, \an \G\in \S_2^s(M) \mapsto \left(g_{\d_1^-(\G)}^\w,e_\G^{(\w,m)}\right) \in G \times E, $$
define a morphism of 2-groupoids; see \ref{ftg}. This result appears in \cite{BS,SW2,FMP1,FMP2}.
\end{Remark}

\subsection{Three-Dimensional holonomy based on a 2-crossed module}

Let $(L \ra{\de} E \ra{\d} G, \t,\left\{,\right \})$ be a 2-crossed module, of Lie groups, and let $(\ll \ra{\de} \le \ra{\d} \lg,\t,\left\{,\right \})$ be the associated differential 2-crossed module.   Recall that we have a left action $\t'$ of $E$ on $L$ defined as $e\t'l=l\{\de(l)^{-1},e\}$, {where $e\in E$ and $l \in L$. The differential form of this action is {$v\t'x=-\{\de(x),v\}$,} where $v \in \le$ and $x \in \ll$}.  Together with the boundary map $\de\colon L \to E$, this defines a crossed module. This permits us to reduce the analysis of 3-dimensional holonomy {based on a 2-crossed module} to the analysis of a  2-dimensional holonomy in the path space, based on a crossed module. {Compare with \ref{greensection} and \cite{FMP1,FMP2,SW2,BS}.}

Suppose we are given forms $\w\in \A^1(M,\lg)$, $m\in \A^2(M,\le)$ and $\theta\in \A^3(M,\ll).$ We suppose that $\d(m)=\W=d\w+[\w,\w]$ and $\de(\theta) =\M=d m+\w \wedge^\t m$. 
Note that $\d\de(\theta)=d \W+\w \wedge^\ad \W=0$, {as it should, by the Bianchi} identity. We define the 3-curvature 4-form $\Theta$ of $(\w,m,\theta)$ as $$\Theta= d \theta+\w\wedge^\t \theta-m\wedge^{\{,\}}m,$$
where the 4-form $m \wedge^{\{,\}} m$ is the antisymmetrisation of the contravariant tensor $6\{m,m\}$. {See the Appendix for notation.}

{{Let  $J\colon [0,1]^3\to M$ be a good 3-path, Definition \ref{good}.}
Put $J(t,s,x)=\G_x(t,s)=\g_s^x(t), \fo t,s,x \in [0,1]$. {A plot $\hat{J}\colon [0,1]^2 \to P(M,*,*')$ is defined as $(s,x) \mapsto \g_s^x$; see subsection \ref{forms}. Here $*=\d_1^-(J)$ and $*'=\d_1^+(J)$.}}

 We define $l_J^{(\w,m,\theta)}(x_0,x)=l_J(x_0,x)$ as being the solution of the differential equation: 
\begin{multline}\label{LJ}
\dx l_J(x_0,x)\\=-l_J(x_0,x)\int_0^1 e_{\G_x}^{(\w,m)}(0,s)\t'{\hat{J}}^*\left(\oint_\w \theta-\oint_\w m *^{\left \{,\right \}} m\right)\left(\ds,\dx\right)ds
\end{multline}
with initial condition $l_J(x_0,x_0)=1_L$.
Compare with (\ref{comp2}).

\subsubsection{The behaviour of the 3-dimensional holonomy on a smooth family of cubes.}
{Let $*,*' \in M$.} Consider the following ($\le$-valued and $\ll$-valued) forms  in the path space $P(M,*,*')$, of smooth paths in $M$ that start in $*$ and finish in $*'$:
$$A=\oint_\w m$$
$$B=-\oint_\w \theta+\oint_\w m *^{\left\{,\right \}}m ,$$
and note, see above,  that $\de(B)$ coincides with the curvature $dA+[A,A]=-\oint_\w \M+\oint_\w m*^{<,>}m$ of $A$. 

Let us calculate the 2-curvature 3-form {$C=dB+A\wedge^{\t'} B$} (in the {path} space) of the pair $(A,B)$. Recall that $(\de\colon \ll \to \le, \t')$ is a differential crossed module. First of all note that, by using the results of \ref{Useful} and the Appendix, it follows that {(taking into account that  $X \t \{u,v\}=\{X \t u,v\}+\{u,X \t v\})$):}
\begin{align*}
d {\oint_\w} m*^{\{,\}}m&=-\oint_\w D_\w m *^{\{,\}} m+\oint_\w m*^{\{,\}} D_\w m-\oint_\w m \wedge^{\{,\}} m\\
&\quad -\oint_\w \W *^\t \left (m *^{\{,\}}m\right) {+\oint_\w m*^{\{,\}} (\W*^\t m)}
\end{align*}
{Thus (by using \ref{Useful} again:)}
\begin{multline*}
C=\oint_\w D_\w \theta
+\oint_\w \W*^\t\theta
-\oint_\w\M*^{\left\{,\right \}}m 
+\oint_\w m*^{\left\{,\right \}}\M 
-\oint_\w m \wedge^{\{,\}} m\\
-\oint_\w \W*^\t \left( m *^{\left\{,\right \}}m\right)+\oint_\w m*^{\{,\}} \left( \W *^{\t}m\right)\\-\oint_\w m \wedge^{\t'}  \oint_\w \theta +\oint_\w m \wedge^{\t'} \oint_\w m *^{\left\{,\right \}}m. \end{multline*}
{Recall that $\M=D_\w m$.
{Then} note that, {since $v \t' x=-\{\de(x),v\}$ for each $v \in \le$ and $x \in \ll$:}}
$$ \oint_\w m \wedge^{\t'}  \oint_\w \theta=-\oint_\w m \wedge^{\{,\}^\op} \oint_\w\de(\theta)= {-\oint_\w m \wedge^{\{,\}^\op} \oint_\w\M}=-\oint \M \wedge^{\{,\}} \oint _\w m.$$
Therefore:
$$ \oint_\w\M*^{\left\{,\right \}}m+\oint_\w m \wedge^{\t'}  \oint_\w \theta=-\oint_\w m*^{\left\{,\right \}^{\op}}\M.$$
{By using condition 6 of the definition of a differential 2-crossed module together with $\de(\theta)=\M$ and $\d(m)=\W$:}
$$\oint_\w \W*^\t\theta+\oint_\w m*^{\left\{,\right \}}\M+\oint_\w m*^{\left\{,\right \}^{\op}}\M=0. $$
Therefore:
\begin{multline*}
C=\oint_\w D_\w \theta -\oint_\w m \wedge^{\{,\}} m
-\oint_\w \W*^\t \left( m *^{\left\{,\right \}}m\right)+\oint_\w m*^{\{,\}} \left( \W *^{\t}m\right)\\+\oint_\w m \wedge^{\t'} \oint_\w m *^{\left\{,\right \}}m. \end{multline*}
{By using  condition 5. of the definition of a differential  2-crossed module follows:}
\begin{align*}
 \oint_\w m \wedge^{\t'} \oint_\w m *^{\left\{,\right \}}m&\doteq-\left( \oint_\w m *^{\left<,\right >}m\right) \wedge^{\{,\}}\oint_\w   m\\
&=-\oint_\w m*^{\{,\}^{\rm op}}\left(m*^{\left<,\right>} m\right)-\oint_\w m*^{\{,\}}\left (m*^{[,]}m\right).
\end{align*}
Thus, {since $\langle u,v\rangle\doteq [u,v]-\d(u) \t v$ for each $u,v \in \le$}:
\begin{multline*}
C=\oint_\w D_\w \theta -\oint_\w m \wedge^{\{,\}} m
-\oint_\w \W*^\t \left( m *^{\left\{,\right \}}m\right)-\oint_\w m*^{\{,\}^{\rm op}} \left( m *^{<,>}m\right) \\ -\oint_{\w}  m *^{\{,\}}\left(m*^{<,>} m\right) .\end{multline*}

By using now  condition 6. of the definition of a 2-crossed module
 it follows that {the} 2-curvature 3-form $C=dA+A\wedge^{\t'} B$ of the pair $(A,B)$ is
$$C=\oint_\w D_\w \theta -\oint_\w m \wedge^{\{,\}} m\doteq \oint_\w \Theta,
$$
recall that $\Theta=D_\w \theta-m\wedge^{\{,\}}m=d \theta+\w\wedge^\t \theta-m\wedge^{\{,\}}m$ denotes the 3-curvature 4-form of $(\w,m,\theta)$. {By using {Corollary} \ref{HOL2} for the crossed module $(\de \colon L \to E,\t')$ (or \cite{FMP1,FMP2}) follows:}
\begin{Theorem} Let $M$ be a smooth manifold and $*,*' \in M$.
Let {$(L\ra{\de} E \ra{\d} G, \t, \{,\}\}$ be a 2-crossed module, and let $(\ll\ra{\de} \le \ra{\d} \lg, \t, \{,\})$} be the associated differential 2-crossed module. Choose forms $\w \in \A^1(M,\lg)$, $m \in \A^2(M,\le)$ and $\theta \in \A^3(M,\ll)$  such that $\de(\theta)=\M=D_\w m$ and $\d(m)=\W$, the curvature of $\w$.

 Let $W\colon D^4 \to M$ be a smooth map such that $\d_1^+(W)=*'$ and $\d_1^-(W)=*${, defining therefore a plot $\hat{W}\colon [0,1]^3 \to P(M,*,*')$.} Let {$W(t,s,x,u)=J_u(t,s,x)={J^x(t,s,u)}=\G_x^u(t,s)$,} where $t,s,x,u \in [0,1]$. {Suppose also that $J_u$ is good  (definition \ref{good}) for all $u$.}  We have:
\begin{multline*}
\du l^{(\w,m,\theta)}_{J_u}(0,1)\\= \int_{[0,1]^3} l^{(\w,m,\theta)}_{J_u}(0,x)\left (e_{\G_x^u}^{(\w,m)}(0,s)\t' {\hat{W}}^* \left ( \oint_\w \Theta \right)\left(\ds,\dx,\du\right)\right) l^{(\w,m,\theta)}_{J_u}(x,1)dsdx\\
-  l^{(\w,m,\theta)}_{J_u}(0,1)\int_0^1 e_{\G_{1}^u}^{(\w,m)}(0,s)\t'{\hat{J^{1}}}^*\left(\oint_\w \theta-\int_\w m *^{\left \{,\right \}} m\right)\left(\ds,\du\right)ds\\
+\left(\int_0^1 e_{\G_{0}^u}^{(\w,m)}(0,s)\t'{\hat{J^0}}^*\left(\oint_\w \theta-{\oint_\w} m *^{\left \{,\right \}} m\right)\left(\ds,\du\right)ds\right)  l^{(\w,m,\theta)}_{J_u}(0,1).
\end{multline*}
\end{Theorem}
\begin{Corollary}\label{3Thin}
 Three dimensional holonomy based on a 2-crossed module is invariant under rank-3 holonomy, restricting to laminated rank-2 holonomy in the boundary. More precisely if $J$ and $J'$ are {good} 3-paths which are rank-3 homotopic (with laminated boundary) then $$l_J^{(\w,m,\theta)}(0,1)=l_{J'}^{(\w,m,\theta)}(0,1),$$
as long as $\de(\theta)=\M=D_\w m$ and $\d(m)=\W$, the curvature of $\w$.
\end{Corollary}

{From now on we will usually abbreviate $l_{J}^{(\w,m,\theta)}(0,1)=l_{J}^{(\w,m,\theta)}$, as we did for $ g_\g^\w\doteq g_\g^\w(0,1)$ and $e_\G^{(\w,m)}\doteq e_\G^{(\w,m)}(0,1)$. }
\subsubsection{The holonomy of the interchange 3-cells}
Fix a 2-crossed module  $(L\ra{\de} E \ra{\d} G, \t, \{,\}\}$ with associated differential  2-crossed module $(\ll\ra{\de} \le \ra{\d} \lg, \t, \{,\})$. As before, consider differential forms  $\w \in \A^1(M,\lg)$, $m \in \A^2(M,\le)$ and $\theta \in \A^3(M,\ll)$  such that $\de(\theta)=\M$ and $\d(m)=\W$, where as usual $\W=d\w+[\w,\w]=d\w+\frac{1}{2}\w \wedge^\ad \w$ and $\M=dm+\w \wedge^\t m$.

Let $\G$ and $\G'$ be two 2-paths with $\d^+_1(\G)=\d^-_1(\G)$. {Let $*=\d^-_1(\G)$, $*'=\d^+_1(\G)$ and $*''=\d^+_1(\G')$.}
\begin{Theorem}[The holonomy of the interchange 3-cell]
We have
$$l^{(\w,m,\theta)}_{\G \#\G'}= e_\G^{(\w,m)} \t' \left \{\left(e_\G^{(\w,m)}\right)^{-1},g^\w_{\d^-_2(\G)}\t e_{\G'}^{(\w,m)} \right\}^{-1}.$$
\end{Theorem}
\begin{Proof}
 Let $J=\G\# \G'$ be the interchange 3-path. Let $\g=\d^-_2(\G)$.  Let $F(s)$ and $F'(s)$ be the paths $\g_s$ and $\g'_s$, respectively, for each $s \in [0,1]${; we have therefore plots $F\colon [0,1] \to P(M,*,*')$ and $F'\colon [0,1] \to P(M,*',*'')$}. As usual {we put} $\g_s(t)=\G(t,s)$ for each $s,t \in [0,1]$, and the same for $\g'_s(t)$. Put $e(x)=e_\G^{(\w,m)}(0,x)$, $e=e(1)$ and $f(x)=e_{\G'}^{(\w,m)}(0,x)$. Also put $l(x)=l_J^{(\w,m,\theta)}(0,x)$. As usual $g_{\g_s}=g_{\g_s}^\w$. Recall $g_{\g_s}=\d(e(s))^{-1}g_\g.$

 By a straightforward explicit calculation, using equation (\ref{LJ}), and the explicit form of $J=\G\#\G'$, indicated by figure \ref{ash}, we have:
\begin{align*}
& \frac{d}{dx} l(x)\\&=l(x)\int_0^1 (g_\g \t f(x)) e(s)\t' 
\left\{  F^*\left({\oint_\w m}  \right)\left(\ds\right) , g_{\g_s}\t F'^*\left(\oint_\w m\right)\left(\dx\right)\right\}ds\\
&=l(x)(g_\g \t f(x)) \t' \int_0^1 e(s)\t' 
\left\{  F^*\left({\oint_\w m}  \right)\left(\ds\right) , g_{\g_s}\t F'^*\left(\oint_\w m\right)\left(\dx\right)\right\}ds.
\end{align*}
We have used the fact that the interchange 3-cell {has} derivative of rank $\leq 2$. Note that by equation (\ref{imppp}):
\begin{multline*}
\int_0^1 e(s)\t' 
\left\{  F^*\left({\oint_\w m}  \right)\left(\ds\right) , g_{\g_s}\t F'^*\left(\oint_\w m\right)\left(\dx\right)\right\}ds\\=
\int_0^1 
\left\{  e(s)F^*\left({\oint_\w m}  \right)\left(\ds\right) , g_{\g_s}\t F'^*\left(\oint_\w m\right)\left(\dx\right)\right\}ds\\
-\int_0^1 
\left\{  e(s), \d\left(F^*\left({\oint_\w m}  \right)\left(\ds\right)\right) g_{\g_s}\t F'^*\left(\oint_\w m\right)\left(\dx\right)\right\}ds.
\end{multline*}
 Since, by definition: $$\frac{d}{ds} {e}(s)={e}(s){F^*\left(\oint_\w m  \right)}\left(\ds\right) \an \frac{d}{ds} {e}^{-1}(s)=-{F^*\left(\oint_\w m  \right)}\left(\ds\right)e^{-1}(s), $$
and also   $\d(e^{-1}(s))g_\g=g_{\g_s}$ (Lemma \ref{Green}), thus
$$\frac{d}{ds} {g}_{\g_s}=-\d\left({F^*\left(\oint_\w m  \right)}\right)g_{\g_s}, $$
it follows, by the Leibnitz rule, and Lemma \ref{trivial}:
\begin{align*}\int_0^1 e(s)\t' 
\Big\{  F^*\left({\oint_\w m}  \right)\left(\ds\right) &, g_{\g_s}\t F'^*\left(\oint_\w m\right)\left(\dx\right)\Big\}ds\\&
=\int_0^1 \frac{\d}{\d s}\left \{e(s),  g_{\g_s}\t F'^*\left(\oint_\w m\right)\left(\dx\right) \right\}
\\&=\left\{e, g_{\g_1}\t F'^*\left(\oint_\w m\right)\left(\dx\right)\right \}.
\end{align*}
{Therefore:}
\begin{equation}\label{AAAAA}
 \frac{d}{dx} l(x)=l(x)(g_\g \t f(x)) \t'\left\{e, g_{\g_1}\t F'^*\left(\oint_\w m\right)\left(\dx\right)\right \}.
\end{equation}

On the other hand, by equation (\ref{Porter}) and the definition of $f(x)$ we have:
\begin{align*}
\frac{d}{dx} &\left\{  e^{-1},  g_\g \t f(x)\right\}= \left\{  e^{-1},  g_\g \t f(x) g_\g \t   F'^*\left(\oint_\w m\right)\left(\dx\right)\right\}\\&=
 \left(e^{-1} \big(g_\g \t f(x)\big) e \right) \t' \left\{e^{-1} ,  g_\g \t F'^*\left(\oint_\w m\right)\left(\dx\right)\right\}\left\{  e^{-1}, g_\g \t f(x)\right\},
\end{align*}
thus in particular:
\begin{align*}
 &\frac{d}{dx}e\t'\left\{  e^{-1}, g_\g \t f(x)\right\}^{-1}\\&=-\left(e\t'\left\{  e^{-1}, g_\g \t f(x)\right\}^{-1}\right)\left(g_\g \t f(x)\right)\t'e \t' \left\{e^{-1}, g_\g \t    F'^*\left(\oint_\w m\right)\left(\dx\right)\right\}\\
&=e \t'\left\{  e, g_{\g} \t f(x)\right\}^{-1}\left(g_\g \t f(x)\right)\t' \left\{e, \d(e^{-1})g_\g \t    F'^*\left(\oint_\w m\right)\left(\dx\right)\right\}\\
&=e \t'\left\{  e, g_{\g} \t f(x)\right\}^{-1}\left(g_\g \t f(x)\right)\t' \left\{e, g_{\g_1} \t    F'^*\left(\oint_\w m\right)\left(\dx\right)\right\},
\end{align*}
taking into account equation (\ref{imppp}) and Lemma \ref{trivial} for the second  step, and  {Lemma} \ref{Green} for the third. By comparing with equation (\ref{AAAAA}) and applying the unicity theorem for ordinary differential equations finishes the proof of the theorem, by using Lemma \ref{trivial} again.
\end{Proof}

\subsubsection{The {Non-abelian} Stokes {Theorem} and the end of the proof of Theorem \ref{Main}}\label{US}

Fix a 2-crossed module $(L\ra{\de} E \ra{\d} G, \t, \{,\}\}$ with associated differential  2-crossed module $(\ll\ra{\de} \le \ra{\d} \lg, \t, \{,\})$. 

Similarly to the  {Non-abelian} Green {Theorem}, the following follows directly from the unicity theorem for ordinary differential equations:
\begin{Lemma}[Non-abelian Stokes Theorem]\label{Stokes}
Let $J$ be a {good} 3-path in $M$. Consider differential forms  $\w \in \A^1(M,\lg)$, $m \in \A^2(M,\le)$ and $\theta \in \A^3(M,\ll)$  such that $\de(\theta)=\M$ and $\d(m)=\W$. 
 We then have:
$$ \d\left(l^{(\w,m,\theta)}_J(x,x')\right)^{-1}e^{(\w,m)}_{\G_x}=
e^{(\w,m)}_{\G_{x'}}$$
recall that we put $J(t,s,x)=\G_x(t,s),\fo t,s,x \in [0,1]$. In addition $e^{(\w,m)}_{\G_x}\doteq e^{(\w,m)}_{\G_x}(0,1)$.
\end{Lemma}

By combining all that was done in this section, it thus follows that the assignment
$$\g\in \S_1(M)\mapsto g_\g^\w \in G,$$ 
$$\G \in \S_2^l(M)\mapsto \left(g^\w_{\d^-_2(\G)}, e_\G^{(\w,m)}\right)$$
and 
$$ J \in \S_3(M) \mapsto \left(g^\w_{\d^-_2(J)}, e_{\d^-_3(J)}^{(\w,m)}, l_J^{(\w,m,\theta)}\right)$$
defines a strict Gray 3-groupoid  map (definition \ref{GrayFunctor}) $$\stackrel{(\w,m,\theta)}{\H}\colon \Sc_3(M) \to \C(\Hc)  ,$$
{where $\C(\Hc)$ is the  Gray 3-groupoid constructed from $\Gc$ (see \ref{nnn}) and $\Sc_3(M)$ is the fundamental Gray 3-groupoid of $M$ (see \ref{thingray})}. What is left to prove is entirely similar to the proof of {Theorem} 39 of \cite{FMP1}; {see \ref{greensection}.}

\section{{An application: Wilson 3-sphere observables}}

\subsection{Complete Whiskering}\label{CWs}
Let $M$ be a manifold. 
Let $J$ be a 3-path with $J(\d D^3)=\{*'\}$. Given a 1-path $\g$ with
 $\d^-_1(\g)=*$  and $\d^+_1(\g)=*'$, let $\g.J$ be the 3-path obtained by filling  $\{z \in D^3, |z|\leq 1/2\}$ with $J$ and the rest of the cube $D^3$ with  $\g$, in the obvious way. This corresponds to the standard way of defining the action of the fundamental groupoid of $M$ on the homotopy  groups $\pi_3(M,*')$. Note that $\g.J$ is well defined up to rank-3 homotopy, with laminated boundary.

Let  $(L\ra{\de} E \ra{\d} G, \t, \{,\}\}$ be a 2-crossed module with associated differential  2-crossed module $(\ll\ra{\de} \le \ra{\d} \lg, \t, \{,\})$. Consider differential forms  $\w \in \A^1(M,\lg)$, $m \in \A^2(M,\le)$ and $\theta \in \A^3(M,\ll)$, as usual  such that $\de(\theta)=\M$ and $\d(m)=\W$.

\begin{Lemma}\label{CW} We have:
$$l_{\g.J}^{(\w,m,\theta)}=g_\g^\w \t l_{J}^{(\w,m,\theta)}.$$
\end{Lemma}
\begin{Proof}
 Follows by the definition of $l_J^{(\w,m,\theta)}$, noting that what was added to $J$ has derivative of rank $\leq 1$, together with the identity $$g\t (e\t' l)=(g\t e) \t' (g \t l),\wh g \in G, e \in E \an l \in L,$$
valid in any 2-crossed module; see equation (\ref{wso}).
\end{Proof}

\subsection{The definition of Wilson 3-sphere observables}

 Let  $S\subset M$ be an oriented 3-sphere $S^3$ embedded in $M$.
 Consider an orientation preserving parametrisation $J \colon D^3/ \d D^3=S^3 \to S \subset M$ of $S$. Define
$$\Wc(S,\w,m,\theta)=l^{(\w,m,\theta)}_J \in \ker \de \subset L,$$
(recall the {Non-abelian} Stokes Theorem), called the Wilson 3-sphere observable.

We have:
\begin{Theorem}
The Wilson 3-sphere observable $\Wc$ does not depend on the parametrisation $J$ of $S$ chosen up to acting by elements of $G$. If $S^*$ denotes $S$ with the reversed orientation we have $\Wc(S^*,\w,m,\theta)=\left(\Wc(S,\w,m,\theta)\right)^{-1}.$
\end{Theorem}
\begin{Proof}
 The second statement is immediate from Theorem \ref{Main}.  Let us prove the first.
Let $J,J'\colon D^3 / \d D^3$ be orientation preserving parametrisations of $S$. Let $J(\d D^3)=*$ and $J'(\d D^3)=*'$. Consider an  isotopy $W\colon D^3 \times I \to S$ connecting $J$ and $J'$, recall that the oriented mapping class group of $S^3$ is trivial. Let $\g(x)=W(\d D^3,x)$, a smooth path in $M$. Obviously {(by using $W$)} $J$ is rank-3 homotopic, with laminated boundary, to $\g. J'$; see subsection \ref{CWs} for this notation. The result follows from {Lemma} \ref{CW} together with {Corollary} \ref{3Thin}.
\end{Proof}

\section{Appendix: Iterated integrals and Chen forms}

Fix a manifold $M$. Given a vector space ${U}$, we denote the vector space of ${U}$-valued  differential $n$-forms in $M$ as $\A^n(M,U)$. Let $V$ and $W$ be vector spaces. Suppose we have a bilinear map $ B\colon (u,v) \in U\times V \mapsto u*v \in W$. If we are given ${U}$ and $V$ valued forms $\a\in \A^a (M,U)$ and $\b \in \A^b(M,V)$ we define  the $W$-valued $(a+b)$-form  $\a \wedge^B\b$ in $M$  (also denoted $\a\wedge^* \b$) as:
$$\a \wedge^B \b= \frac{(a+b)!}{a!b!}{\rm Alt}(\a \otimes^B \b) \in \A^{a+b}(M,W).$$
Here $\a\otimes^B \b$ is the covariant tensor $B\circ (\a \times \b)$ and ${\rm Alt}$ denotes the natural projection from the vector space of  $W$-valued  covariant tensor fields in $M$ onto the vector space of $W$-valued differential forms in $M$. Note that if $B^{\op}$ denotes the bilinear map $V \times U \to W$ such that $B^{\op}(v,u)=B(u,v)$ then 
$$\b\wedge^{B^\op} \a=(-1)^{ab} \a \wedge^B \a.$$
On the other hand if we define  {$J(\w)={(-1)}^m \w$,} where $\w$ is an $m$-form, we will continue to have $d(\a \wedge\b)=(d \a )\wedge^B \b+J(\a) \wedge^B d \b$.

\subsection{Iterated integrals}\label{iterated}
{For details see \cite{Ch,BT}. Let ${\U}$  be an open set in some $\R^n$. Let ${\w}_1,\ldots ,{\w}_n$ be forms in  $\A^{n_i+1}( \R \times \U,W_i)$, {for $i=1, \ldots , k$,} where $W_i$ are vector spaces.  Suppose we are given bilinear maps $ W_i \times W_{i+1} \to W_{i+1}$, say $(v,w) \mapsto v*w$.}
Given $-\infty<a<b<+\infty$, we define 
$\a=\oint_a^b {\w}_1dt=\oint_a^b \w_1 $ as being the $W_1$-valued $n_1$-form in ${\U}$ such that
$${\a(v_1,\ldots,v_{n_1})=\int_a^b {\w}_1\left(\dt, v_1,\ldots,v_{n_1}\right)dt .}$$
Having defined the $W_m$-valued $(n_1+\ldots+n_m)$-form $\oint_{a}^b {\w}_1*\ldots *{\w}_{m} dt$, then the  $W_{m+1}$-valued $(n_1+\ldots+n_m+n_{m+1})$-form  $\oint_a^b {\w}_1*\ldots *{\w}_{m+1} dt$ is $$\oint_a^b {\w}_1*\ldots *{\w}_{m+1}=\int_a^b \left(\oint_a^t  {\w}_1*\ldots *{\w}_{m} dt'\right) \wedge^* \i_{\dt}({\w}_{m+1}) dt,$$
where $\i_X(\w)$ denotes the contraction of a form $\w$ with a vector field $X$. Sometimes parentheses may be inserted to denote the order in which we apply bilinear maps (if not the order above).

For simplicity (but with enough generality for this article), suppose that all forms ${\w}_i$ vanish when $t=0$ and $t=1$. Given  $b \in [0,1]$, let $\w_b=r_b^*(\w)$, where $r_b(x)=(b,x)$ for each $x \in \U$.  As in \cite{Ch} we can prove that:
\begin{align}
d \oint_0^b \w&=-\int_0^b (\i_{\dt} d\w)dt+\int_0^b \frac{d}{dt} \w_t dt
             =\w_b-\w_0 -\int_0^b (\i_{\dt} d\w) dt\\\label{Ch1}
             &=\w_b-\oint_0^b (d  \w)dt.
\end{align}
{In general in the case when $V_i=\R$ we have (by induction):}
\begin{multline}
d \oint_0^b {\w}_1*\ldots *{\w}_{m} dt =\sum_{i=1}^m (-1)^i \oint_0^bJ{\w}_1*\ldots* J {\w}_{i-1} *(d{\w}_i) *{\w}_{i+1}*\ldots* {\w}_m
\\+\sum_{i=1}^{m-1} (-1)^i \oint_0^bJ{\w}_1*\ldots* J {\w}_{i-1} *( J{\w}_i  \wedge^*{\w}_{i+1})*{\w}_{i+2}*\ldots* {\w}_m\\-(-1)^m\left(\oint_0^b  J({\w}_1)*\ldots*J({\w}_{m-1})\right) \wedge^* r_b^*(\w_m)
\end{multline}
For a proof see \cite{Ch}.

From now on we use the following convention 
$$\oint \w_1*\ldots *\w_n=\oint_0^1 \w_1*\ldots *\w_n dt.$$

\subsection{Forms in the space of curves}\label{forms}
{Let $M$ be a smooth manifold. Let $P(M,*,*')$ denote the space of all smooth curves $\g\colon [0,1] \to M$  {that start in $*$ and finish at $*'$}. Recall that a plot is a map $F\colon {\U \to} P(M,*,*')$, such that the associated map $F'\colon I \times \U \to M$ {given by}  $F'(t,x)=F(x)(t)$ is smooth; see \cite{Ch}. (Here $\U$ is an open set in some $\R^n$.)}

 By definition a $p$-form $\a$ in {$P(M,*,*')$} is given by a rule which associates a $p$-form $F^*(\a)$ in ${\U}$ to each plot {$F\colon \U \to P(M,*,*')$,}  satisfying the following compatibility condition:
for any smooth map $g\colon \U' \to \U$, where $\U'$ is an open set in some $\R^n$, we have $(F\circ g)^*(\w)=g^*(F^*(\w))$, in $\U'$. The sum, exterior product and exterior derivative of forms in the space of curves are defined {as:
$$F^*(\a+\b)=F^*(\a)+F^*(\b),\, \, F^*(\a\wedge \b)=F^*(\a)\wedge F^*(\b), \,\, F^*(d \a)=d F^*(\a),$$
respectively, }
for each plot {$F\colon \U \to P(M,*,*')$.}

Let ${\w}_1,\ldots,{\w}_n$ be forms in $M$ of degrees $a_i+1$, {for $i=1,\ldots , n$.} Then we define a $a_1+\ldots+a_n$-form $\oint {\w}_1*\ldots *{\w}_n$ in the path space $P(M,*,*')$ by putting:{
$$F^*\left(\oint {\w}_1*\ldots* {\w}_n\right)=\oint (F')^*({\w}_1)*\ldots *(F')^*({\w}_n).$$ 
}

{Note that if we are given a smooth map $f\colon [0,1]\times [0,1]^n \to M$, such that $f(\{0\} \times [0,1]^n)=*$ and $f(\{1\} \times [0,1]^n)=*'$, then we have  an associated plot $\hat{f} \colon [0,1]^n \to P(M,*,*')$, sometimes denoted by $F \colon [0,1]^n \to P(M,*,*')$.}

\section*{Acknowledgements}
The first author  was   supported by the Centro de Matem\'{a}tica da
Universidade do Porto {\it www.fc.up.pt/cmup}, financed by {\em Funda\c{c}\~{a}o para a Ci\^{e}ncia e a Tecnologia} (FCT)  through the programmes POCTI and POSI, with Portuguese and European Community structural funds. This work was  partially supported by the {\em Programa Operacional  Ci\^{e}ncia e Inova\c{c}\~{a}o 2010},  financed by FCT and 
cofinanced by the European Community fund FEDER, in part through the research project Quantum Topology POCI/MAT/60352/2004.

{We would like to thank   Tim Porter, Marco Zambom, Urs Schreiber, Branislav Jur\v{c}o and Bj\"{o}rn Gohla  for useful discussions  and / or  comments.}

\end{document}